\documentclass[11pt,a4paper]{article}

\usepackage{amsmath,amssymb,amsthm}
\usepackage{mathtools}
\usepackage[font=small, labelfont=bf]{caption}
\usepackage[labelformat=simple]{subcaption}

\usepackage{fullpage}
\usepackage[colorlinks=true, citecolor=blue]{hyperref}
\usepackage{tikz}
\usetikzlibrary{knots}
\usetikzlibrary{backgrounds}
\usepackage{authblk}

\RequirePackage[sc]{mathpazo}
\RequirePackage[T1]{fontenc}

\newcommand{\tube}{\mathbb T}
\newcommand{\tubeLM}{\mathbb T_{L,M}}
\newcommand{\pham}{p^\mathrm{H}}
\newcommand{\kappaham}{\kappa^\mathrm{H}}

\DeclareMathOperator{\trunk}{trunk}

\newtheorem{result}{Result}

\author[1]{N R Beaton\thanks{Corresponding author, \href{mailto:nrbeaton@unimelb.edu.au}{nrbeaton@unimelb.edu.au}}}
\author[2]{J W Eng}
\author[3]{K Ishihara}
\author[4]{K Shimokawa}
\author[2]{C E Soteros}
\affil[1]{School of Mathematics and Statistics, The University of Melbourne, Australia}
\affil[2]{Department of Mathematics and Statistics, University of Saskatchewan, Saskatoon, Canada}
\affil[3]{Faculty of Education, Yamaguchi University, Yamaguchi, Japan}
\affil[4]{Department of Mathematics, Saitama University, Saitama, Japan}

\title{Characterising knotting properties of polymers in nanochannels}

\begin{document}
\maketitle

\begin{abstract}
Using a lattice model of polymers in a tube, we define one way to characterise different configurations of a given knot as either ``local'' or ``non-local'' and, for several ring polymer models, we provide both theoretical and numerical evidence that, at equilibrium, the non-local configurations are more likely than the local ones. These characterisations are based on a standard approach for measuring the ``size'' of a knot within a knotted polymer chain. The method involves associating knot-types to subarcs of the chain, and then identifying a knotted subarc with minimal arclength; this arclength is then the knot-size. If the resulting knot-size is small relative to the whole length of the chain, then the knot is considered to be localised or ``local''. If on the other hand the knot-size is comparable to the length of the chain, then the knot is considered to be ``non-local''.

Using this definition, we establish that all but exponentially few sufficiently long self-avoiding polygons (closed chains)  in a tubular sublattice of the simple cubic lattice are ``non-locally'' knotted. This is shown to also hold for the case when the same polygons are subject to an external tensile force, as well as in the extreme case when they are as compact as possible (no empty lattice sites). We also provide numerical evidence for small tube sizes that at equilibrium non-local knotting is more likely than local knotting, regardless of the strength of the stretching or compressing force. We note however that because of the tube confinement, the occurrence of non-local knotting  in walks (open chains)  is significantly different than for polygons. The relevance of these results to recent experiments involving DNA knots in solid-state nanopores is also discussed.
\end{abstract}

\section{Introduction}

Motivated in part by experimental studies of  DNA packing  in viral capsids~\cite{Arsuaga_2005,Marenduzzo_2009}  and DNA translocation through nanopores~\cite{Dai_2016,Micheletti_2012,Plesa_2016,Suma_2017}, there has been much recent interest in understanding and characterising the entanglement complexity of confined polymers, and determining any dependencies on the extent or the mechanism of confinement.  For knots in polymers, one measure of interest 
has been the  average ``size'' of the knotted part of the polymer.    With such a measurement, one can then characterise the knotting as ``local''  when the size of the knotted part is small compared to the whole length of the polymer, or otherwise  as ``non-local''.   Based on polymer scaling theory and supporting numerical evidence, it is generally accepted that  local knotting is dominant for unconfined polymers, as polymer length grows, while computer simulation  studies of knotting  in collapsed or spherically confined polymers suggest that non-local knotting dominates~\cite{Micheletti_2011}.  Recent simulation and experimental studies~\cite{Micheletti_2012, orlandini_statics_2018, Plesa_2016, Suma_2017} of an intermediate regime of confinement, where the polymer is confined to a channel, tube or pore, such that polymer growth or motion is tightly restricted in  two spatial dimensions but unrestricted in the third, has indicated that local knotting dominates,  unlike when all three spatial dimensions are restricted. Some differences between the entanglement complexity of  open versus closed chains were also observed.  

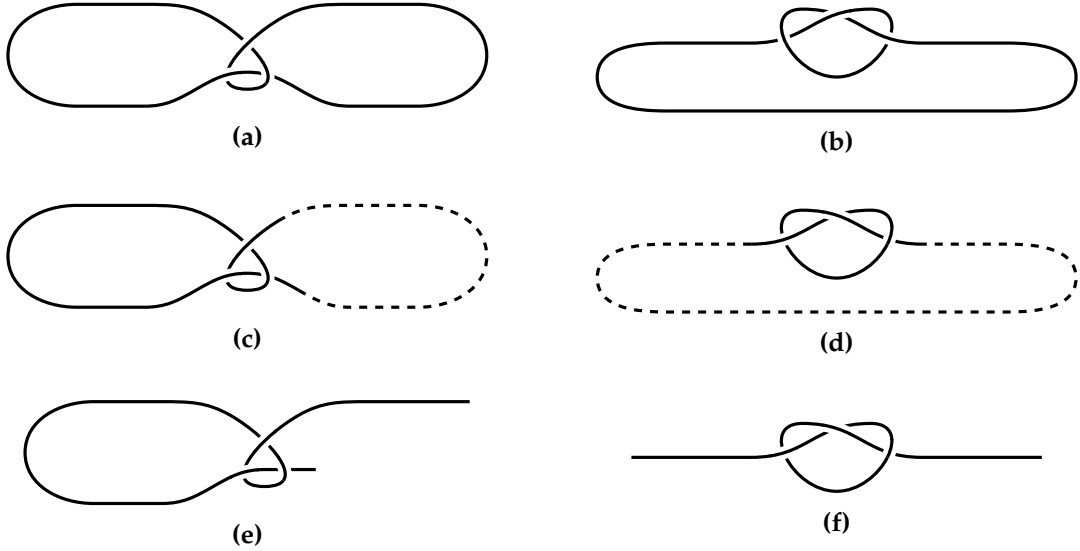
\begin{figure}
\centering
\begin{subfigure}{0.48\textwidth}
\centering
\begin{tikzpicture}[scale=0.45]
\begin{knot}[consider self intersections=true, clip width=4, clip radius=5pt, end tolerance=1pt]
\strand[very thick, black] (1.5,-0.5) to [out=-30, in=180] (3,-1) to [out=0, in=180] (5,-1) to [out=0, in=-90] (7,0.5) to [out=90, in=0] (5,2) to [out=180, in=0] (3,2) to [out=180, in=30] (1,1.6) to [out=-150, in=180, looseness=1.5] (0,-0.5) to [out=0, in=-30, looseness=1.5] (-1,1.6) to [out=150, in=0] (-3,2) to [out=180, in=0] (-5,2) to [out=180, in=90] (-7,0.5) to [out=-90, in=180] (-5,-1) to [out=0, in=180] (-3,-1) to [out=0, in=-150] (-1.5,-0.5) to [out=30, in=180] (0,0) to [out=0, in=150] (1.5,-0.5);
\flipcrossings{2}
\end{knot}
\end{tikzpicture}
\caption{ }
\label{fig:nonlocal_loop}
\end{subfigure}
\begin{subfigure}{0.48\textwidth}
\centering
\begin{tikzpicture}[scale=0.45]
\begin{knot}[consider self intersections=true, clip width=4, clip radius=5pt, end tolerance=0.5pt]
\strand[very thick, black] (0,0) to [out=0, in=0, looseness=1.5] 
(1,2) to [out=180, in=0, looseness=1.2] 
(-2.5,1) to [out=180, in=0] 
(-5,1)  to [out=180, in=90] 
(-7,0) to [out=-90, in=180] 
(-5,-1) to [out=0, in=180] 
(5,-1) to [out=0, in=-90] 
(7,0) to [out=90, in=0] 
(5,1) to [out=180, in=0] 
(2.5,1) to [out=180, in=0, looseness=1.2] 
(-1,2) to [out=180, in=180, looseness=1.5] (0,0);
\flipcrossings{1,3}
\end{knot}
\end{tikzpicture}
\caption{ }
\label{fig:local_loop}
\end{subfigure}

\vspace{0.5cm}
\begin{subfigure}{0.48\textwidth}
\centering
\begin{tikzpicture}[scale=0.45]
\begin{knot}[consider self intersections=true, clip width=4, clip radius=5pt, end tolerance=1pt]
\strand[very thick, black] (1,1.6) to [out=-150, in=180, looseness=1.5] (0,-0.5) to [out=0, in=-30, looseness=1.5] (-1,1.6) to [out=150, in=0] (-3,2) to [out=180, in=0] (-5,2) to [out=180, in=90] (-7,0.5) to [out=-90, in=180] (-5,-1) to [out=0, in=180] (-3,-1) to [out=0, in=-150] (-1.5,-0.5) to [out=30, in=180] (0,0) to [out=0, in=150] (1.5,-0.5);
\strand[very thick, black, dashed] (1.5,-0.5) to [out=-30, in=180] (3,-1) to [out=0, in=180] (5,-1) to [out=0, in=-90] (7,0.5) to [out=90, in=0] (5,2) to [out=180, in=0] (3,2) to [out=180, in=30] (1,1.6);
\flipcrossings{2}
\end{knot}
\end{tikzpicture}
\caption{ }
\label{fig:nonlocal_loop_colored}
\end{subfigure}
\begin{subfigure}{0.48\textwidth}
\centering
\begin{tikzpicture}[scale=0.45]
\begin{knot}[consider self intersections=true, clip width=4, clip radius=5pt, end tolerance=1pt]
\strand[very thick, black] (2.5,1) to [out=180, in=0, looseness=1.2] (-1,2) to [out=180, in=180, looseness=1.5] (0,0) to [out=0, in=0, looseness=1.5] (1,2) to [out=180, in=0, looseness=1.2] (-2.5,1);
\strand[very thick, black, dashed] (-2.5,1) to [out=180, in=0] (-5,1) to [out=180, in=90] (-7,0) to [out=-90, in=180] (-5,-1) to [out=0, in=180] (5,-1) to [out=0, in=-90] (7,0) to [out=90, in=0] (5,1) to [out=180, in=0] (2.5,1);
\flipcrossings{1,3}
\end{knot}
\end{tikzpicture}
\caption{ }
\label{fig:local_loop_colored}
\end{subfigure}

\vspace{0.5cm}
\begin{subfigure}{0.48\textwidth}
\centering
\begin{tikzpicture}[scale=0.45]
\begin{knot}[consider self intersections=true, clip width=4, clip radius=5pt, end tolerance=1pt]
\strand[very thick, black] (6,2) to [out=180, in=0] (3,2) to [out=180, in=30] (1,1.6) to [out=-150, in=180, looseness=1.5] (0,-0.5) to [out=0, in=-30, looseness=1.5] (-1,1.6) to [out=150, in=0] (-3,2) to [out=180, in=0] (-5,2) to [out=180, in=90] (-7,0.5) to [out=-90, in=180] (-5,-1) to [out=0, in=180] (-3,-1) to [out=0, in=-150] (-1.5,-0.5) to [out=30, in=180] (0,0) to [out=0, in=180] (1.5,0);
\flipcrossings{2}
\end{knot}
\end{tikzpicture}
\caption{ }
\label{fig:nonlocal_loop_walk}
\end{subfigure}
\begin{subfigure}{0.48\textwidth}
\centering
\begin{tikzpicture}[scale=0.45]
\begin{knot}[consider self intersections=true, clip width=4, clip radius=5pt, end tolerance=0.5pt]
\strand[very thick, black] (6,1) to [out=180, in=0] (2.5,1) to [out=180, in=0, looseness=1.2] (-1,2) to [out=180, in=180, looseness=1.5] (0,0) to [out=0, in=0, looseness=1.5] (1,2) to [out=180, in=0, looseness=1.2] (-2.5,1) to [out=180, in=0] (-6,1);
\flipcrossings{1,3}
\end{knot}
\end{tikzpicture}
\caption{ }
\label{fig:local_loop_walk}
\end{subfigure}
\caption{Illustrations of \subref{fig:nonlocal_loop} a non-local trefoil and \subref{fig:local_loop} a local trefoil, with knotted subarcs drawn with solid lines in \subref{fig:nonlocal_loop_colored} and \subref{fig:local_loop_colored}; open chains with \subref{fig:nonlocal_loop_walk} a non-local trefoil component and \subref{fig:local_loop_walk} a local trefoil component.}
\label{figstrings}
\end{figure}

In this paper we explore this latter type of confinement further using lattice models of both closed (self-avoiding polygons) and open (self-avoiding walk) chains in tubular subsets of the simple cubic lattice.   One major advantage of the lattice model approach is that we are able to prove results related to knot localisation for arbitrary tube dimensions and, for the case of small tube sizes, we are able to perform exact calculations related to this.   We study an equilibrium model of polymers in a tube subject to a tensile force $f$ and prove results about the limiting free energy and the likelihood of occurrence of different knotted patterns as a function of polymer length. Through this we provide both theoretical and numerical evidence that for a knotted ring polymer at equilibrium in a nanochannel or nanopore,  a knot configuration such as that shown in Figure~\ref{fig:nonlocal_loop} is more likely than that shown in Figure~\ref{fig:local_loop}, regardless of the strength of the force $f$ and whether or not it is a compressing ($f<0$) or a stretching ($f>0$) force.   Furthermore,  knot configurations such as that in Figure~\ref{fig:nonlocal_loop} are on average tighter than those in Figure~\ref{fig:local_loop} and hence would be expected to translocate through a nanopore faster.   
We also establish for the lattice model that, as observed experimentally, the situation is quite different for open chains, with knot configurations like that of Figure~\ref{fig:nonlocal_loop_walk} being rare (in fact exponentially rare) compared to that shown in Figure~\ref{fig:local_loop_walk}, because of the significant entropic disadvantage associated with the formation of a long bend.   

To obtain and explain these results, it is necessary to first characterise the differences between the different knotted configurations in Figure \ref{figstrings}.  We do this here by using two different measures of knot size.   For example, when one measures the size of the knotted part according to the size of the region in which the crossings are concentrated, then the configurations shown in both
Figures~\ref{fig:nonlocal_loop} and~\ref{fig:local_loop} correspond to examples of  ``local'' or ``tight'' knotting and this is consistent, for example, with a measure of how long it would take the knotted part to translocate through a nanopore.  However, using another standard measure for determining knot-size, namely using the length of a smallest knotted subarc,  leads to characterising Figure~\ref{fig:nonlocal_loop} as ``non-local'' knotting (a knotted subarc is drawn with a solid line in Figure~\ref{fig:nonlocal_loop_colored}); while that in Figure~\ref{fig:local_loop} is ``local'' (a knotted subarc is drawn with a solid line in Figure~\ref{fig:local_loop_colored}).    Using this latter measure, our results indicate that non-local knotting is more likely than local knotting in a tube.  This suggests that the likelihood of non-local knotting in other restricted spaces, like spheres, will be dependent on the definition chosen for knot-size.

The remainder of this paper is organised as follows. In Section~\ref{sec:themodels} we present definitions of the models under consideration.  In Section~\ref{sec:characterising}, definitions for  two measures of knot size and our classification scheme for knotted patterns are given.  After that, we present exact and numerical results about the models. Then we present the known theoretical results (outlines of proofs are given in the appendix (Section~\ref{ssec:proofs})) and the methods used.   
Finally in Section~\ref{sec:conclusion} we have a detailed discussion of the results and their relevance to other numerical and experimental studies.

\section{The models}\label{sec:themodels}

We will use a general model for polygons in lattice tubes subject to an external force which has been studied previously;  the notation and definitions used here (unless stated otherwise) are as in~\cite{Beaton_2016}.

{For non-negative integers $L,M$, let} \(\tubeLM\equiv\tube\subset \mathbb Z^3\) be the semi-infinite $L\times M$ tube on the cubic lattice defined by
\[\tube = \{(x,y,z)\in \mathbb Z^3:x\geq0, 0\leq y\leq L, 0 \leq z \leq M\}.\]
Define $\mathcal P_\tube$ to be the set of self-avoiding polygons in $\tube$ which occupy at least one vertex in the plane $x=0$, and let $\mathcal P_{\tube,n}$ be the subset of $\mathcal P_\tube$ comprising polygons with $n$ edges ($n$ even). Then let $p_{\tube,n}  = |\mathcal P_{\tube,n}|$. 

\begin{figure}
\centering
\begin{subfigure}{0.45\textwidth}
\resizebox{\textwidth}{!}{
\begin{tikzpicture}[rotate around x=270, scale=1.7]
\tikzset{vblue/.style={circle, draw=blue, line width=2.5pt, fill=blue, inner sep=2pt}}
\tikzset{vblack/.style={circle, draw=black, line width=2.5pt, fill=black, inner sep=2pt}}
\tikzset{vred/.style={circle, draw=red, line width=2.5pt, fill=red, inner sep=2pt}}
\node[vblack] at (0,0,0) {}; \node[vblack] at (0,0,1) {}; \node[vblack] at (0,1,1) {}; \node[vblack] at (0,2,1) {};
\node[vblue] at (1,0,0) {}; \node[vblue] at (1,1,0) {}; \node[vblue] at (1,2,0) {}; \node[vblue] at (1,0,1) {}; \node[vblue] at (1,1,1) {}; \node[vblue] at (1,2,1) {};
\node[vblack] at (2,0,0) {}; \node[vblack] at (2,1,0) {}; \node[vblack] at (2,2,0) {}; \node[vblack] at (2,0,1) {}; \node[vblack] at (2,2,1) {};
\node[vblack] at (3,0,0) {}; \node[vblack] at (3,1,0) {}; \node[vblack] at (3,2,0) {}; \node[vblack] at (3,0,1) {}; \node[vblack] at (3,2,1) {};
\node[vblack] at (4,0,0) {}; \node[vblack] at (4,1,0) {}; \node[vblack] at (4,2,0) {}; \node[vblack] at (4,0,1) {}; \node[vblack] at (4,1,1) {}; \node[vblack] at (4,2,1) {};
\node[vblack] at (5,0,0) {}; \node[vblack] at (5,1,0) {}; \node[vblack] at (5,2,0) {}; \node[vblack] at (5,0,1) {}; \node[vblack] at (5,1,1) {}; \node[vblack] at (5,2,1) {};
\node[vblack] at (6,2,0) {}; \node[vblack] at (6,0,0) {}; \node[vblack] at (6,1,0) {};
\node[vblack] at (6,0,1) {};
\begin{scope}[on background layer]
\begin{knot}[consider self intersections=true, clip width=3, clip radius=4pt, end tolerance=1pt]
\strand[line width=2pt, black] (2,1,0) -- (2,0,0) -- (1,0,0) -- (0,0,0) -- (0,0,1) -- (0,1,1) -- (0,2,1) -- (1,2,1) -- (2,2,1);
\strand[line width=2pt, black] (1,0,1) -- (2,0,1);
\strand[line width=2pt, black] (1,2,0) -- (2,2,0);
\strand[line width=2pt, blue] (1,0,1) -- (1,1,1) -- (1,1,0) -- (1,2,0);
\strand[line width=2pt, green] (2,1,0) -- (3,1,0);
\strand[line width=2pt, green] (2,2,0) -- (3,2,0);
\strand[line width=2pt, green] (2,0,1) -- (3,0,1);
\strand[line width=2pt, green] (2,2,1) -- (3,2,1);
\strand[line width=2pt, black] (3,0,1) -- (3,0,0) -- (4,0,0) -- (4.5,0,0);
\strand[line width=2pt, black] (3,2,1) -- (4,2,1) -- (4.5,2,1);
\strand[line width=2pt, black] (3,2,0) -- (4,2,0) -- (4.5,2,0);
\strand[line width=2pt, black] (3,1,0) -- (4,1,0) -- (4,1,1) -- (4,0,1) -- (4.5,0,1);
\strand[line width=2pt, black] (4.5,0,0) -- (5,0,0) -- (5,1,0) -- (5,1,1) -- (5,2,1) -- (4.5,2,1);
\strand[line width=2pt, black] (4.5,0,1) -- (5,0,1) -- (5.5,0,1);
\strand[line width=2pt, black] (4.5,2,0) -- (5,2,0) -- (5.5,2,0);
\strand[line width=2pt, black] (5.5,0,1) -- (6,0,1) -- (6,0,0) -- (6,1,0) -- (6,2,0) -- (5.5,2,0);
\flipcrossings{2,3,4,5,6,8}
\end{knot}
\end{scope}
\node at (7,2,1) {};
\end{tikzpicture}
}
\caption{}
\label{fig:polygon_2x1_hinge_section_block}
\end{subfigure}

\vspace{0.25cm}
\begin{subfigure}
{0.45\textwidth}
\resizebox{\textwidth}{!}{
\begin{tikzpicture}[rotate around x=270, scale=1.7]
\tikzset{vblack/.style={circle, draw=black, line width=2.5pt, fill=black, inner sep=2pt}}
\tikzset{vblack/.style={circle, draw=black, line width=2.5pt, fill=black, inner sep=2pt}}
\tikzset{vblack/.style={circle, draw=black, line width=2.5pt, fill=black, inner sep=2pt}}
\node[vblack] at (0,0,0) {}; \node[vblack] at (0,0,1) {}; \node[vblack] at (0,1,1) {}; \node[vblack] at (0,2,1) {};
\node[vblack] at (1,0,0) {}; \node[vblack] at (1,1,0) {}; \node[vblack] at (1,2,0) {}; \node[vblack] at (1,0,1) {}; \node[vblack] at (1,1,1) {}; \node[vblack] at (1,2,1) {};
\node[vblack] at (2,0,0) {}; \node[vblack] at (2,1,0) {}; \node[vblack] at (2,2,0) {}; \node[vblack] at (2,0,1) {}; \node[vblack] at (2,2,1) {};
\node[vblack] at (3,0,0) {}; \node[vblack] at (3,1,0) {}; \node[vblack] at (3,2,0) {}; \node[vblack] at (3,0,1) {}; \node[vblack] at (3,2,1) {};
\node[vblack] at (4,0,0) {}; \node[vblack] at (4,1,0) {}; \node[vblack] at (4,2,0) {}; \node[vblack] at (4,0,1) {}; \node[vblack] at (4,1,1) {}; \node[vblack] at (4,2,1) {};
\node[vblack] at (5,0,0) {}; \node[vblack] at (5,1,0) {}; \node[vblack] at (5,2,0) {}; \node[vblack] at (5,0,1) {}; \node[vblack] at (5,1,1) {}; \node[vblack] at (5,2,1) {};
\node[vblack] at (6,2,0) {}; \node[vblack] at (6,0,0) {}; \node[vblack] at (6,1,0) {};
\node[vblack] at (6,0,1) {};
\begin{scope}[on background layer]
\begin{knot}[consider self intersections=true, clip width=3, clip radius=5.5pt, end tolerance=1pt]
\strand[line width=2pt, black] (0,0,0) -- (1,0,0) -- (2,0,0) -- (2,1,0) -- (3,1,0) -- (4,1,0) -- (4,1,1) -- (4,0,1) -- (5,0,1) -- (6,0,1) -- (6,0,0) -- (6,1,0) -- (6,2,0) -- (5,2,0) -- (4,2,0) -- (3,2,0) -- (2,2,0) -- (1,2,0) -- (1,1,0) -- (1,1,1) -- (1,0,1) -- (2,0,1) -- (3,0,1) -- (3,0,0) -- (4,0,0) -- (5,0,0) -- (5,1,0) -- (5,1,1) -- (5,2,1) -- (4,2,1) -- (3,2,1) -- (2,2,1) -- (1,2,1) -- (0,2,1) -- (0,1,1) -- (0,0,1) -- cycle;
\strand[line width=2pt, red] (0.5,0,-0.15) -- (0.5,0,0.15);
\strand[line width=2pt, red] (0.5,2,0.85) -- (0.5,2,1.15);
\strand[line width=2pt, red] (5.5,0,0.85) -- (5.5,0,1.15);
\strand[line width=2pt, red] (5.5,2,-0.15) -- (5.5,2,0.15);
\flipcrossings{1,2,3,6,8,9,10,11,12}
\end{knot}
\end{scope}
\node at (7,2,1) {};
\end{tikzpicture}
}
\caption{}
\label{fig:polygon_2x1_2strings}
\end{subfigure}

\vspace{0.25cm}
\begin{subfigure}
{0.45\textwidth}
\resizebox{\textwidth}{!}{
\begin{tikzpicture}[rotate around x=270, scale=1.7]
\tikzset{vblack/.style={circle, draw=black, line width=2.5pt, fill=black, inner sep=2pt}}
\tikzset{vblack/.style={circle, draw=black, line width=2.5pt, fill=black, inner sep=2pt}}
\tikzset{vblack/.style={circle, draw=black, line width=2.5pt, fill=black, inner sep=2pt}}
\node[vblack] at (0,0,0) {}; \node[vblack] at (0,0,1) {}; \node[vblack] at (0,1,1) {}; \node[vblack] at (0,2,1) {};
\node[vblack] at (1,0,0) {}; \node[vblack] at (1,1,0) {}; \node[vblack] at (1,2,0) {}; \node[vblack] at (1,0,1) {}; \node[vblack] at (1,1,1) {}; \node[vblack] at (1,2,1) {};
\node[vblack] at (2,0,0) {}; \node[vblack] at (2,1,0) {}; \node[vblack] at (2,2,0) {}; \node[vblack] at (2,0,1) {}; \node[vblack] at (2,2,1) {};
\node[vblack] at (3,0,0) {}; \node[vblack] at (3,1,0) {}; \node[vblack] at (3,2,0) {}; \node[vblack] at (3,0,1) {}; \node[vblack] at (3,2,1) {};
\node[vblack] at (4,0,0) {}; \node[vblack] at (4,1,0) {}; \node[vblack] at (4,2,0) {}; \node[vblack] at (4,0,1) {}; \node[vblack] at (4,1,1) {}; \node[vblack] at (4,2,1) {};
\node[vblack] at (5,0,0) {}; \node[vblack] at (5,1,0) {}; \node[vblack] at (5,2,0) {}; \node[vblack] at (5,0,1) {}; \node[vblack] at (5,1,1) {}; \node[vblack] at (5,2,1) {};
\node[vblack] at (6,2,0) {}; \node[vblack] at (6,0,1) {}; \node[vblack] at (6,1,1) {}; \node[vblack] at (6,2,1) {};  \node[vblack] at (6,0,0) {}; \node[vblack] at (6,1,0) {};
\node[vblack] at (7,1,1) {}; \node[vblack] at (7,2,1) {};
\begin{scope}[on background layer]
\begin{knot}[consider self intersections=true, clip width=3, clip radius=5.5pt, end tolerance=1pt]
\strand[line width=2pt, black] (0,0,0) -- (1,0,0) -- (2,0,0) -- (2,1,0) -- (3,1,0) -- (4,1,0) -- (4,1,1) -- (4,0,1) -- (5,0,1) -- (6,0,1) -- (6,0,0) -- (6,1,0) -- (6,2,0) -- (5,2,0) -- (4,2,0) -- (3,2,0) -- (2,2,0) -- (1,2,0) -- (1,1,0) -- (1,1,1) -- (1,0,1) -- (2,0,1) -- (3,0,1) -- (3,0,0) -- (4,0,0) -- (5,0,0) -- (5,1,0) -- (5,1,1) -- (6,1,1) -- (7,1,1) -- (7,2,1) -- (6,2,1) -- (5,2,1) -- (4,2,1) -- (3,2,1) -- (2,2,1) -- (1,2,1) -- (0,2,1) -- (0,1,1) -- (0,0,1) -- cycle;
\strand[line width=2pt, red] (0.5,0,-0.15) -- (0.5,0,0.15);
\strand[line width=2pt, red] (0.5,2,0.85) -- (0.5,2,1.15);
\strand[line width=2pt, red] (6.5,1,0.85) -- (6.5,1,1.15);
\strand[line width=2pt, red] (6.5,2,0.85) -- (6.5,2,1.15);
\flipcrossings{1,2,3,7,8,9,10,11,12}
\end{knot}
\end{scope}
\end{tikzpicture}
}
\caption{}
\label{fig:polygon_2x1_2stringslocal}
\end{subfigure}

\caption{\subref{fig:polygon_2x1_hinge_section_block}~A 36-edge polygon $\pi$ that fits inside $\tube_{L,M}$ with $L\geq 2$ and $M\geq 1$;  the tube extends without bound to the right and the span $s(\pi)=6$. Blue vertices and edges denote the hinge $H_1$ of $\pi$, and green edges denote the section $S_3$ of $\pi$. \subref{fig:polygon_2x1_2strings}~The locations of the two pairs of vertical red lines indicate the locations of the two 2-sections in this polygon; in this example, the polygon can be decomposed into a start unknot pattern, a proper trefoil knot pattern, and an end unknot pattern. The proper knot pattern is classified as non-local in this case. \subref{fig:polygon_2x1_2stringslocal}~A local proper knot pattern in the same tube with span 7.}
\label{fig2tubedef}
\end{figure}
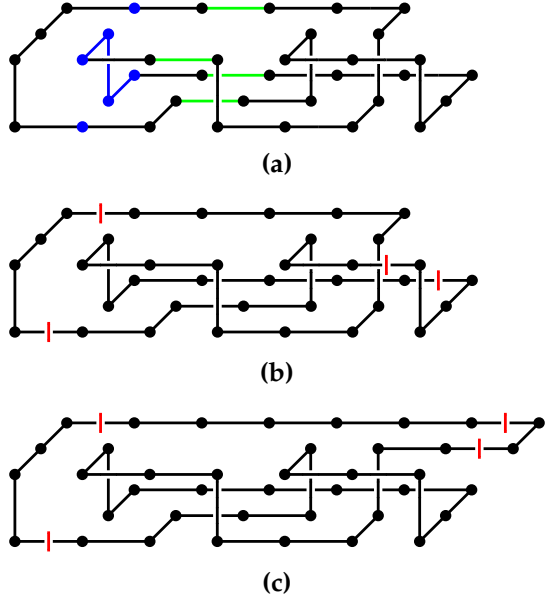

We define the \emph{span} $s(\pi)$ of a polygon $\pi\in\mathcal P_\tube$ to be the maximal $x$-coordinate reached by any of its vertices and we use $|\pi|$ to denote the number of edges in $\pi$. 
See Figure~\ref{fig:polygon_2x1_hinge_section_block} for a polygon $\pi$ that fits in a $2\times1$ tube with $s(\pi)=6$ and $|\pi|=36$.
To model a force acting parallel to the $x$-axis, we  associate a fugacity (Boltzmann weight) $e^{fs(\pi)}$ with each polygon $\pi$. Let $p_{\tube,n}(s)$ be the number of polygons in $\mathcal P_{\tube,n}$ with span $s$. Then the ``fixed-edge'' model partition function is given by
\begin{displaymath}
Z_{\tube,n}(f)  = \sum_{|\pi| = n} e^{fs(\pi)} = \sum_s p_{\tube,n}(s) e^{fs}.
\end{displaymath}
Thus $f\ll 0$ corresponds to the ``compressed'' regime while $f\gg 0$ corresponds to the ``stretched'' regime.
For this model, the probability of a polygon $\pi \in \mathcal P_{\tube,n}$ is given by
\begin{displaymath}
{\mathbb{P}}_{n}^{({\rm{ed}}, f)}(\pi)=\frac{e^{fs(\pi)}}{Z_{\tube,n}(f)  }  .
\end{displaymath}

The (limiting) \emph{free energy per edge} of polygons in $\tube$ is defined as
\[\mathcal F_\tube(f) = \lim_{n\to\infty}\frac1n \log Z_{\tube,n}(f).\]
This is known~\cite{Atapour_2009} to exist for all $f$.

For $f=0$, it has been proved that~\cite{Soteros_1988, Soteros_1989}
\begin{align}
\mathcal F_\tube(0) &= \lim_{n\to\infty} n^{-1}\log p_{\tube,n} \notag \\
&< \lim_{n\to\infty} n^{-1}\log c_{\tube,n} \notag \\
&< \lim_{n\to\infty} n^{-1}\log p_n = \lim_{n\to\infty} n^{-1}\log c_n \equiv \kappa \equiv \log \mu, \label{eqnmudef}
\end{align}
where $c_n$ is the number of $n$-step self-avoiding walks (SAWs) in $\mathbb Z^3$ starting at the origin and $\kappa$ is their connective constant, and $c_{\tube,n}$ is the number of these confined to $\tube$.

A subset of self-avoiding polygons in $\tube$ are \emph{Hamiltonian} polygons: those which occupy every vertex in a $s\times L\times M$ subtube of $\tube$. These serve as an idealised model of tightly packed ring polymers, in addition to being a useful lower bound for general polygons in the $f<0$ compressed regime. We define the number of Hamiltonian polygons, $\pham_{\tube,n}$, to be the number of  $n$-edge polygons in $\mathcal{P}_{\tube,n}$ which have span $s$ and occupy every vertex in an $s\times L\times M$ subtube of $\tube$. 
We define $W=(L+1)(M+1)$ (the number of vertices in an integer plane $x=i\geq 0$ of the tube)  and will assume without loss of generality that $L\geq M$;
note that $\pham_{\tube,n} = 0$ if $n$ is not a multiple of $W$.
The following limit has been proved to exist~\cite{Beaton_2016} (see also~\cite{Eng2014MSc}):
\begin{equation*}\label{eqn:ham_growthrate1}
\kappaham_\tube \equiv \lim_{s\to\infty}\frac{1}{(s+1)W} \log \pham_{\tube,(s+1)W}.
\end{equation*}
Furthermore, using this, $\mathcal F_\tube(f) $, the free energy per edge, is bounded as follows:
\begin{align}
\max\{{f/2},(f/W)+ \kappaham_\tube\} &\leq \mathcal F_\tube(f) \leq  \max\{{f/W},{f/2}\} +\mathcal F_\tube(0),
\label{eqnfreeenergybounds}
\end{align}
with $\mathcal F_\tube(f)$ asymptotic to the lower bound 
for  $f\to\infty$ for any $\tube$, and for $f\to -\infty$ for small tube sizes (this is conjectured to be true for any $\tube$), see~\cite{Beaton_2016}.

Here, we will also be interested in the dual  model, called the ``fixed-span'' model, with partition function given by
\[Q_{\tube,s}(g) = \sum_n p_{\tube,n}(s) e^{g n}.\]
For this partition function, when $g \gg 0$  densely packed (in terms of number of edges per span) polygons dominate the partition function, while when $g \ll 0$ polygons with very few edges per span dominate.  
For this model, the probability associated with a span $s$ polygon $\pi$ is given by
\begin{displaymath}
{\mathbb{P}}_{s}^{({\rm{sp}}, g)}(\pi)=\frac{e^{g|\pi|}}{Q_{\tube,s}(g)  }  
\end{displaymath}
and the associated (limiting) {\it free energy per span} exists~\cite{Beaton_2016} (see also~\cite{Atapour2008PhD}):
\[\mathcal G_\tube(g) = \lim_{s\to\infty}\frac1s \log Q_{\tube,s}(g).\]

Both models correspond to special cases of the grand canonical partition function 
\[G_{\tube}(f,g) = \sum_s \sum_n p_{\tube,n}(s) e^{g n+ f s}\]
and can be studied using transfer matrix methods~\cite{Atapour2008PhD}.  

Hamiltonian polygons can also be studied using transfer matrices~\cite{Beaton_2016,Eng2014MSc} and we will also be interested in the fixed-span model where polygons are restricted to being Hamiltonian and the probability associated with a span $s$ Hamiltonian polygon $\pi$ is given by
\begin{displaymath}
{\mathbb{P}}_{s}^{(\mathrm{H})}(\pi)=\frac{1}{\pham_{\tube,(s+1)W} }  .
\end{displaymath}

\section{Characterising local knotting and classifying knotted patterns}\label{sec:characterising}

Given a polygon of prime knot-type $K$, one standard approach for measuring the ``size'' of the knotted part in the polygon is to find a minimal length sub-walk of the polygon which has knot type $K$, and then define the size of the knot to be the length of this sub-walk. This of course requires a method for assigning a knot-type to an open chain; there are various ways to do so (see~\cite{Micheletti_2011,StuEnzo} for reviews). If $K$ is composite, then one can use a similar approach for each of the components. We refer to this measurement of knot-size as  \emph{arclength knot-size}.

Motivated by the arclength knot-size definition, in this section we introduce an approach for classifying the knotted parts of polygons in $\tube$ as either ``local'' or ``non-local''.
We will also introduce another measure for knot-size which is particularly suited to polygons in narrow tubes. References will be made to Figure~\ref{fig2tubedef} for illustration.

Let $\pi$ be a polygon  in $\tube$ with span $s$. Thus $\pi$ is embedded in $\tube$ between the planes $x=0$ and $x=s$. For half-integers $k\in\mathbb Z+\frac12$ with $0<k<s$, we say that $\pi$ has a \emph{2-section} at $x=k$ if the plane $x=k$ intersects $\pi$ at exactly two points. (Equivalently, $\pi$ has exactly two edges in the $x$-direction between $x=k-\frac12$ and $x=k+\frac12$.) 
See Figure~\ref{fig:polygon_2x1_2strings} where $k=1/2$ or $k=11/2$.
If $\pi$ has $m$ 2-sections, let $t(\pi) = (t_1,\ldots,t_m)$ be the (ordered) set of $x$-values at which they occur. Clearly $m\leq s$; if $m=0$ then $t(\pi)$ is empty.

The 2-sections of a polygon $\pi$ in $\tube$ naturally partition it into a sequence of ``segments''. Moreover, if $\pi$ has prime knot-type, then typically the segment of $\pi$ which contains  the ``knotted part"  will lie between two successive 2-sections. See for examples Figures~\ref{fig:polygon_2x1_2strings}  and \ref{fig:polygon_2x1_2stringslocal} for polygons with trefoil knot-type.  It is this idea which will allow us to locate, measure and classify knot components within polygons in $\tube$.

Take $\pi$ in $\tube$ with $m\geq 2$. The 2-sections of $\pi$ partition it into a sequence of segments that we call \emph{connect-sum patterns} (\emph{cs-patterns} for short).
Then for any $1 \leq i < m$, the segment of $\pi$ between $x=t_i$ and $x=t_{i+1}$ 
%(including the ``half-edges'' at either end) 
is called a  \emph{proper} cs-pattern of $\pi$. If $m\geq1$ then the segment of $\pi$ between $x=0$ and $x=t_1$ is called the \emph{start} cs-pattern of $\pi$; likewise, the segment of $\pi$ between $x=t_m$ and $x=s(\pi)$ is called the \emph{end} cs-pattern of $\pi$. See Figures~\ref{fig:polygon_2x1_2strings}  and \ref{fig:polygon_2x1_2stringslocal} for examples with $m=2$. Note that proper, start and end cs-patterns are examples of, respectively, proper, left-most and right-most patterns as defined in~\cite{Beaton_2016}. Here we define the \emph{span} of a  (resp.~start, end) cs-pattern to be $t_{i+1}-t_i+1$ (resp.~$t_1+\frac12$, $s(\pi)-t_m+\frac12$).

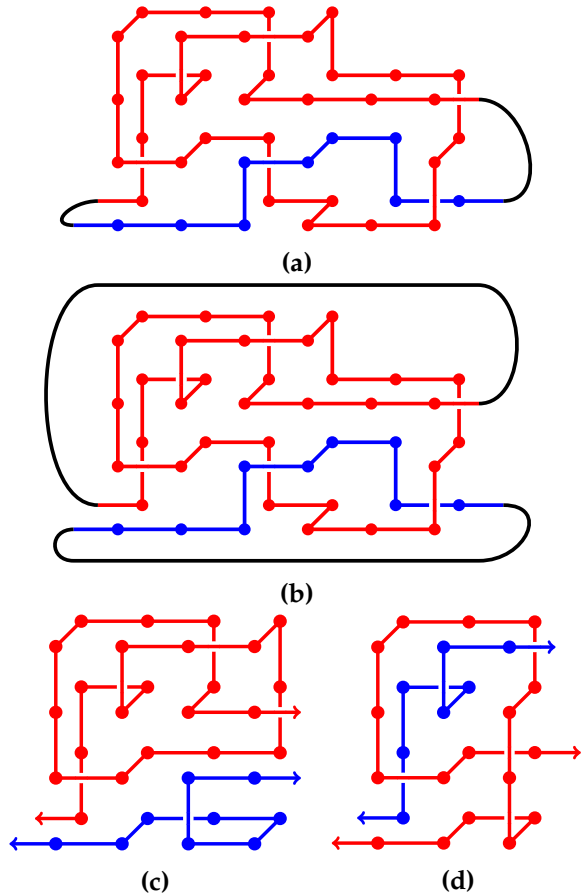
\begin{figure}
\centering
\begin{subfigure}{0.47\textwidth}
\resizebox{\textwidth}{!}{
\begin{tikzpicture}[rotate around x=270, scale=1.5]
\tikzset{vblue/.style={circle, draw, fill=blue, inner sep=2pt}}
\tikzset{vred/.style={circle, draw, fill=red, inner sep=2pt}}
\begin{knot}[consider self intersections=true, clip width=3, clip radius=4pt, end tolerance=1pt]
\strand[line width=2.5pt, blue] (-0.7,0,0) -- (0,0,0) node [vblue] {} -- (1,0,0) node [vblue] {} -- (2,0,0) node [vblue] {} -- (2,0,1) node [vblue] {} -- (3,0,1) node [vblue] {} -- (3,1,1) node [vblue] {} -- (4,1,1) node [vblue] {} -- (4,1,0) node [vblue] {} -- (5,1,0) node [vblue] {} -- (5.7,1,0);
\strand[line width=2.5pt, red] (-0.7,1,0) -- (0,1,0) node [vred] {} -- (0,1,1) node [vred] {} -- (0,1,2) node [vred] {} -- (1,1,2) node [vred] {} -- (1,0,2) node [vred] {} -- (1,0,3) node [vred] {} -- (2,0,3) node [vred] {} -- (3,0,3) node [vred] {} -- (3,1,3) node [vred] {} -- (3,1,2) node [vred] {} -- (4,1,2) node [vred] {}-- (5,1,2) node [vred] {} -- (5,1,1) node [vred] {} -- (5,0,1) node [vred] {} -- (5,0,0) node [vred] {} -- (4,0,0) node [vred] {} -- (3,0,0) node [vred] {} -- (3,1,0) node [vred] {} -- (2,1,0) node [vred] {} -- (2,1,1) node [vred] {} -- (1,1,1) node [vred] {} -- (1,0,1) node [vred] {} -- (0,0,1) node [vred] {} -- (0,0,2) node [vred] {} -- (0,0,3) node [vred] {} -- (0,1,3) node [vred] {} -- (1,1,3) node [vred] {} -- (2,1,3) node [vred] {} -- (2,1,2) node [vred] {} -- (2,0,2) node [vred] {} -- (3,0,2) node [vred] {} -- (4,0,2) node [vred] {} -- (5,0,2) node [vred] {} -- (5.7,0,2);
\flipcrossings{2,3,4,6}
\end{knot}

\draw[line width=2.5pt, black] (-0.7,0,0) to [out=180, in=180, looseness=2] (-0.7,1,0);
\draw[line width=2.5pt, black] (5.7,1,0) to [out=0, in=0, looseness=1.2] (5.7,0,2);

\node at (-1.5,0,0) {};
\node at (6,3,1) {};
\end{tikzpicture}
}
\caption{}
\label{fig:nonlocal-5_1-numerator}
\end{subfigure}

\begin{subfigure}{0.47\textwidth}
\resizebox{\textwidth}{!}{
\begin{tikzpicture}[rotate around x=270, scale=1.5]
\tikzset{vblue/.style={circle, draw, fill=blue, inner sep=2pt}}
\tikzset{vred/.style={circle, draw, fill=red, inner sep=2pt}}
\begin{knot}[consider self intersections=true, clip width=3, clip radius=4pt, end tolerance=1pt]
\strand[line width=2.5pt, blue] (-0.7,0,0) -- (0,0,0) node [vblue] {} -- (1,0,0) node [vblue] {} -- (2,0,0) node [vblue] {} -- (2,0,1) node [vblue] {} -- (3,0,1) node [vblue] {} -- (3,1,1) node [vblue] {} -- (4,1,1) node [vblue] {} -- (4,1,0) node [vblue] {} -- (5,1,0) node [vblue] {} -- (5.7,1,0);
\strand[line width=2.5pt, red] (-0.7,1,0) -- (0,1,0) node [vred] {} -- (0,1,1) node [vred] {} -- (0,1,2) node [vred] {} -- (1,1,2) node [vred] {} -- (1,0,2) node [vred] {} -- (1,0,3) node [vred] {} -- (2,0,3) node [vred] {} -- (3,0,3) node [vred] {} -- (3,1,3) node [vred] {} -- (3,1,2) node [vred] {} -- (4,1,2) node [vred] {}-- (5,1,2) node [vred] {} -- (5,1,1) node [vred] {} -- (5,0,1) node [vred] {} -- (5,0,0) node [vred] {} -- (4,0,0) node [vred] {} -- (3,0,0) node [vred] {} -- (3,1,0) node [vred] {} -- (2,1,0) node [vred] {} -- (2,1,1) node [vred] {} -- (1,1,1) node [vred] {} -- (1,0,1) node [vred] {} -- (0,0,1) node [vred] {} -- (0,0,2) node [vred] {} -- (0,0,3) node [vred] {} -- (0,1,3) node [vred] {} -- (1,1,3) node [vred] {} -- (2,1,3) node [vred] {} -- (2,1,2) node [vred] {} -- (2,0,2) node [vred] {} -- (3,0,2) node [vred] {} -- (4,0,2) node [vred] {} -- (5,0,2) node [vred] {} -- (5.7,0,2);
\flipcrossings{2,3,4,6}
\end{knot}

\draw[line width=2.5pt, black] (-0.7,0,0) to [out=180, in=180, looseness=2] (-0.7,0,-0.5) to [out=0, in=180] (5.7,0,-0.5) to [out=0, in=0, looseness=2] (5.7,1,0);
\draw[line width=2.5pt, black] (-0.7,1,0) to [out=180, in=180, looseness=0.8] (-0.7,1,3.5) to [out=0, in=180] (5.3,1,3.5) to [out=0, in=0, looseness=1.1] (5.7,0,2);

\node at (-1.5,0,0) {};
\node at (6,3,1) {};
\end{tikzpicture}
}
\caption{}
\label{fig:nonlocal-5_1-denominator}
\end{subfigure}

\begin{subfigure}{0.241\textwidth}
\resizebox{\textwidth}{!}
{\begin{tikzpicture}[rotate around x=270, scale=2]
\tikzset{vblue/.style={circle, draw, fill=blue, inner sep=3pt}}
\tikzset{vred/.style={circle, draw, fill=red, inner sep=3pt}}
\begin{knot}[consider self intersections=true, clip width=3, clip radius=5pt, end tolerance=1pt]
\strand[line width=3pt, blue, <->] (-0.7,0,0) -- (0,0,0) node [vblue] {} -- (1,0,0) node [vblue] {} -- (1,1,0) node [vblue] {} -- (2,1,0) node [vblue] {} -- (3,1,0) node [vblue] {} -- (3,0,0) node [vblue] {} -- (2,0,0) node [vblue] {} -- (2,0,1) node [vblue] {} -- (3,0,1) node [vblue] {} -- (3.7,0,1);
\strand[line width=3pt, red, <->] (-0.7,1,0) -- (0,1,0) node [vred] {} -- (0,1,1) node [vred] {} -- (0,1,2) node [vred] {} -- (1,1,2) node [vred] {} -- (1,0,2) node [vred] {} -- (1,0,3) node [vred] {} -- (2,0,3) node [vred] {} -- (3,0,3) node [vred] {} -- (3,1,3) node [vred] {} -- (3,1,2) node [vred] {} -- (3,1,1) node [vred] {} -- (2,1,1) node [vred] {} -- (1,1,1) node [vred] {} -- (1,0,1) node [vred] {} -- (0,0,1) node [vred] {} -- (0,0,2) node [vred] {} -- (0,0,3) node [vred] {} -- (0,1,3) node [vred] {} -- (1,1,3) node [vred] {} -- (2,1,3) node [vred] {} -- (2,1,2) node [vred] {} -- (2,0,2) node [vred] {} -- (3,0,2) node [vred] {} -- (3.7,0,2);
\flipcrossings{1,2,3,5}
\end{knot}
\end{tikzpicture}}
\caption{}
\label{fig:local_trefoil}
\end{subfigure}
\hspace{1ex}
\begin{subfigure}{0.209\textwidth}
\resizebox{\textwidth}{!}
{\begin{tikzpicture}[rotate around x=270, scale=2]
\tikzset{vblue/.style={circle, draw, fill=blue, inner sep=3pt}}
\tikzset{vred/.style={circle, draw, fill=red, inner sep=3pt}}
\begin{knot}[consider self intersections=true, clip width=3, clip radius=5pt, end tolerance=1pt]
\strand[line width=3pt, red, <->] (-0.7,0,0) -- (0,0,0) node [vred] {} -- (1,0,0) node [vred] {} -- (1,1,0) node [vred] {} -- (2,1,0) node [vred] {} --  (2,0,0) node [vred] {} -- (2,0,1) node [vred] {} -- (2,0,2) node [vred] {} -- (2,1,2) node [vred] {} -- (2,1,3) node [vred] {} -- (1,1,3) node [vred] {} -- (0,1,3) node [vred] {} -- (0,0,3) node [vred] {} -- (0,0,2) node [vred] {} -- (0,0,1) node [vred] {} -- (1,0,1) node [vred] {} -- (1,1,1) node [vred] {} -- (2,1,1) node [vred] {} -- (2.7,1,1);
\strand[line width=3pt, blue, <->] (-0.7,1,0) -- (0,1,0) node [vblue] {} -- (0,1,1) node [vblue] {} -- (0,1,2) node [vblue] {} -- (1,1,2) node [vblue] {} -- (1,0,2) node [vblue] {} -- (1,0,3) node [vblue] {} -- (2,0,3) node [vblue] {} -- (2.7,0,3);
\flipcrossings{1,3,5}
\end{knot}
\end{tikzpicture}
}
\caption{}
\label{fig:nonlocal_trefoil}
\end{subfigure}
\caption{\subref{fig:nonlocal-5_1-numerator}~ An illustration of the denominator closure of a proper cs-pattern $\sigma$.  The blue strand corresponds to $\sigma_1$ and the red to $\sigma_2$ and their union is proper cs-pattern $\sigma$.  The denominator closure  is obtained by adding the black arcs and yields a closed curve with knot-type $5_1$, i.e. $DC(\sigma)=5_1$. \subref{fig:nonlocal-5_1-denominator}~The numerator closure of the same pattern. The numerator closure gives a link with one component a $3_1\,\,(\neq 5_1)$ knot and the other an unknot ($0_1\neq 5_1$); hence this is a non-local knot pattern. Here $NC_1(\sigma)=0_1$ and $NC_2(\sigma)=3_1$. \subref{fig:local_trefoil}~A local trefoil knot pattern that can occur in a  Hamiltonian polygon. \subref{fig:nonlocal_trefoil}~A non-local trefoil knot pattern that can occur in a  Hamiltonian polygon. }
\label{fig3_non-localdef}
\end{figure}

Any proper cs-pattern $\sigma$ (between $x=t_i$ and $x=t_{i+1}$, for some $i\geq 1$) is the union of two ``strands'' (self-avoiding walks) $\sigma_1$ and $\sigma_2$, each extending from the left end of the pattern to the right. On the left side of the plane $x=t_i$,  joining the two left ends of $\sigma_1$ and $\sigma_2$ to each other and then, on the right side of the plane $x=t_{i+1}$, joining their two right ends to each other, yields what  we call the \emph{denominator closure} of the cs-pattern. See Figure~\ref{fig:nonlocal-5_1-numerator}. \emph{Note here that if the overall polygon has knot-type $K$ and one of its proper cs-patterns has denominator closure of knot-type $K'$, then $K'$ must be part of the knot-decomposition of $K$.} Let $DC(\sigma)$ be the knot-type of the denominator closure of $\sigma$.  If $DC(\sigma)\neq 0_1$, then we say that $\sigma$ is a \emph{knot pattern} with 
knot-type $DC(\sigma)$.

Alternatively, the two endpoints of $\sigma_1$ (resp.~$\sigma_2$) can be {\it  reconnected} to each other (outside of $\tube$) to form a (possibly separable) link. We call this the \emph{numerator closure} of $\sigma$. See Figure~\ref{fig:nonlocal-5_1-denominator}. Here we are not interested in the overall link-type of the numerator closure; we instead only care about the knot-types of its two components. Let $NC_1(\sigma)$ (resp.~$NC_2(\sigma)$) be the knot-type of the closure of $\sigma_1$ (resp.~$\sigma_2$).

We are now prepared to give our definitions of local and non-local knot patterns. Let $\sigma$ be a cs-pattern of a polygon. If $\sigma$ is a proper pattern with $DC(\sigma) = K\neq 0_1$ but such that 
$NC_1(\sigma)\neq K\# K'$ and $NC_2(\sigma)\neq K\# K''$ for any $K'$ or $K''$ (i.e. $K$ is not in the knot decomposition of either $NC_1(\sigma)$ or  $NC_2(\sigma)$), 
then the knot $K$ cannot be discovered by examining only one of the strands of $\sigma$. In this case both strands are needed to detect $K$ and hence since the two strands of $\sigma$ could potentially be far apart along the contour of the entire polygon, we define $\sigma$ to be a \emph{non-local} knot pattern. For example, the denominator closure in Figure~\ref{fig:nonlocal-5_1-numerator} is a $5_1$ knot while the components of the numerator closure in Figure~\ref{fig:nonlocal-5_1-denominator} are a trefoil and an unknot; the corresponding proper cs-pattern is therefore a non-local knot pattern.   For all other cases we classify $\sigma$ as a \emph{local} knot pattern. Note that there are examples of cs-patterns which we classify here as ``local'' which might be more appropriately classified as ``non-local'' (see Appendix Section~\ref{ssec:comments} for a discussion), however, for small tube sizes such cs-patterns have a relatively low probability of occurrence.

%There are three special cases which do not fit into the above definitions. A proper start pattern can only be closed in one way, and since the arclength of the pattern will (usually) be small relative to the length of the whole polygon, we regard all start patterns which are knots as local knot patterns. The same applies to end patterns.

%The third special case is when a polygon contains no 2-strings at all. Such polygons are in fact exponentially rare in all tube sizes (this is straightforward to prove, but we will not discuss it further here). Nevertheless, for the purposes of this article, we will classify knotted polygons with no 2-strings as non-local knot patterns.

The decomposition of polygons into cs-patterns 
leads to a second measure for knot-size in a polygon as follows.  A polygon with prime knot-type $K$ has at most one knot pattern (either local or non-local). If it has a knot pattern, then we define its {\it connect-sum knot-size} to be the number of edges in the knot pattern.
If there is no knot pattern in the polygon, then the length of the whole polygon is used.   A similar approach can be used for a composite knot $K$, where we use the number of edges in each (if any) knot pattern to determine a total connect-sum knot-size.

In the next section we explore numerically occurrence probabilities associated with local and non-local cs-patterns for small tube sizes and specific knot-types and we provide evidence that non-local knot patterns are more probable than local ones.  
However, even though a polygon $\pi$ with knot-type $K\neq 0_1$ contains a knot pattern of a given type (local or non-local),  this does not necessarily guarantee that the polygon is locally or non-locally knotted (using the standard arclength classification); this aspect is explored more thoroughly and theoretically in Section~\ref{sec:theoretical}.   

\section{Exact and simulation results}

\begin{figure}
\centering
\resizebox{0.5\textwidth}{!}{
%\tikzsetnextfilename{figure4}
%\input{fg_plot_v3.file}
\includegraphics{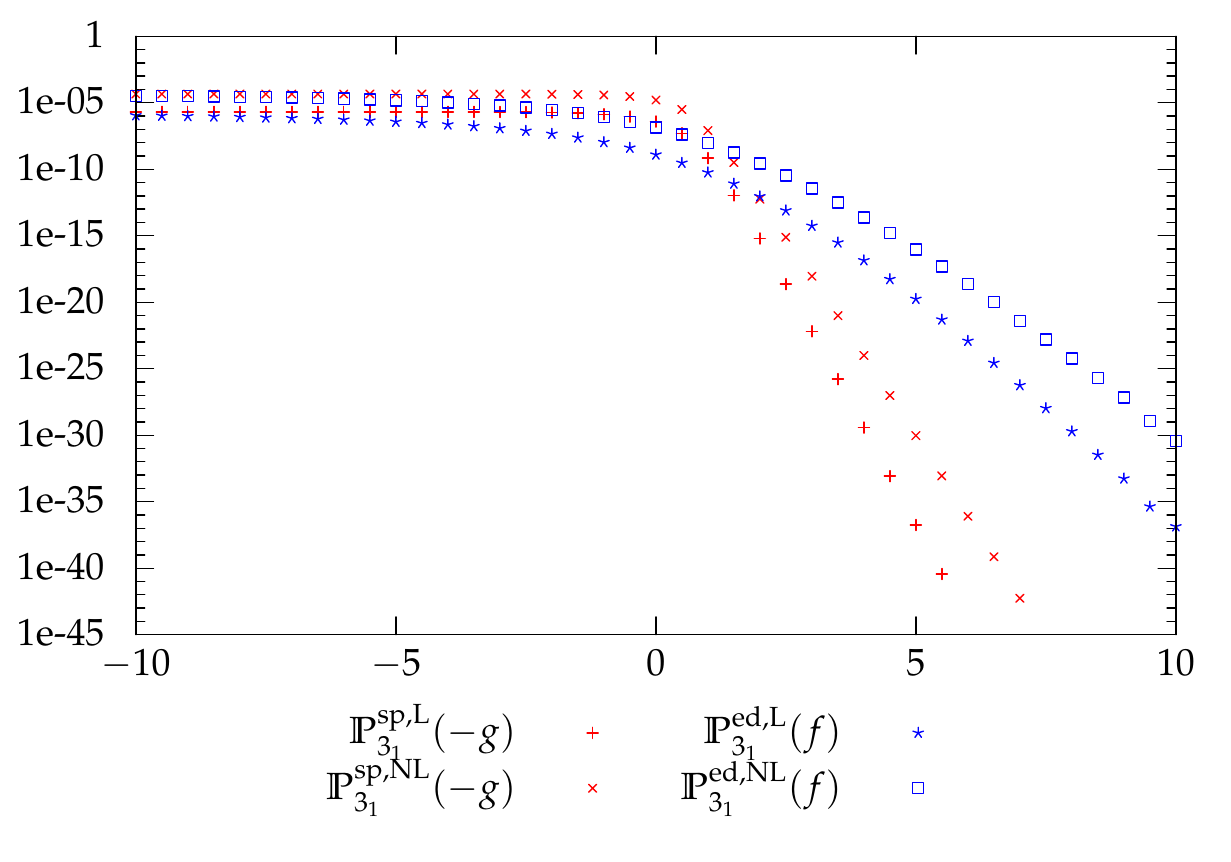}
}
\caption{Log scale plot of the probabilities of the smallest local and non-local trefoil patterns in the $3\times1$ tube, as functions of $f$ (blue) and $-g$ (red).}
\label{fig_fg}
\end{figure}

Exact generation was used to determine all smallest-span non-local and local trefoil knot patterns in the $2\times1$ and $3\times1$ tubes.
(See Figure~\ref{fig2tubedef} for examples of such patterns in the $2\times 1$ tube and Figure~\ref{fig3_non-localdef} in the $3\times 1$ tube.) 
Counts 
are shown in Table~\ref{table:jeremy_data}.     In all cases the number of non-local knot patterns greatly exceeds the number of local knot patterns, suggesting that non-local trefoil knot patterns may be more likely to occur than local ones.   To explore whether this conclusion depends on the model used (fixed-edge or fixed-span), limiting probabilities of occurrence of each type of pattern were determined under each of the distributions ${\mathbb{P}}_{n}^{({\rm{ed}}, f)}$ ($-\infty <f <\infty$, $n\to\infty$) and ${\mathbb{P}}_{s}^{({\rm{sp}}, g)}$ ($-\infty <g<\infty$, $s\to\infty$).  These limiting probabilities can be determined (see for example~\cite{Eng2014MSc}) 
 from the eigenvalues and eigenvectors of the transfer-matrix.  Figure~\ref{fig_fg} shows the results for the $3\times 1$ tube.  In this figure, 
for the fixed-edge model (${\mathbb{P}}_{n}^{({\rm{ed}}, f)}$),  $\mathbb{P}_{3_1}^{\textup{ed},\textup{L}}(f)$ denotes the limiting ($n\to\infty$) probability of occurrence of a smallest local trefoil knot pattern at a section of a polygon and $\mathbb{P}_{3_1}^{\textup{ed},\textup{NL}}(f)$ denotes the corresponding probability for the non-local patterns. 
Similarly, for the fixed-span model (${\mathbb{P}}_{s}^{({\rm{sp}}, g)}$),  $\mathbb{P}_{3_1}^{\textup{sp},\textup{L}}(g)$ denotes the limiting ($s\to\infty$) probability of occurrence of a smallest local trefoil knot pattern at a section of a polygon and $\mathbb{P}_{3_1}^{\textup{sp},\textup{NL}}(g)$ denotes the corresponding probability for the non-local patterns. 
Further note that Figure~\ref{fig_fg} shows the results for the fixed-edge probabilities with the horizontal axis corresponding to  $f$ while for the fixed-span probabilities it corresponds to $-g$.  
The latter was done to make an easier comparison between the models, since positive values of $f$ and negative values of $g$ both have a stretching effect on polygons.
Although not shown here, the observed trends were similar for the $2\times 1$ tube. The results are summarised below.

\begin{result}\label{result:smallest-trefoil-patterns}
For $\tube_{2,1}$ and $\tube_{3,1}$,  the limiting occurrence probability
of the smallest-span non-local trefoil knot patterns is greater than that
for the smallest-span local trefoil knot patterns, regardless of the strength of the force $f$ or the edge density weight $g$.
\end{result}

\begin{table}
\caption{Numbers of trefoil patterns of smallest spans in the $2\times1$ and $3\times1$ tubes, for all and Hamiltonian polygons.}
\centering
\begin{tabular}{crrrrr} \hline
	Tube & & & & Ham. & Ham. \\
	Size & Span & Non-local & Local & Non-local & Local \\
	\hline\hline
	$2\times1$		& 6 & 116 		& 0		& 32		& 0\\
				& 7 & 5,584	& 304	& 668	& 80 \\
				& 8 & 141,292	& 14,932	& 8,020	&  1,388 \\
	\hline
	$3\times1$		& 4 &	1,964	&0		& 232 & 0 \\
				& 5 &	762,984	&29,272	& 17,568 	  & 1,448 \\
	\hline
\end{tabular}
\label{table:jeremy_data}
\end{table}

Determining the knot-type of a polygon or a knot pattern requires the whole polygon or knot pattern.  However, since the numbers of polygons and knot patterns in a tube grow exponentially with either span or the number of edges,  we limited the exact generation analysis to the case of the smallest trefoil patterns.  To explore further whether the trend observed for the smallest trefoil patterns holds for other knots,  a Monte Carlo approach was developed to generate random polygons in the tube, based on a method of~\cite{Alm_1990}. The Monte Carlo method is also based on transfer-matrices and can be used to generate a set of independent and identically distributed polygons from any of the distributions
$\{{\mathbb{P}}_{n}^{({\rm{ed}}, f)},{\mathbb{P}}_{s}^{({\rm{sp}}, g)}\}$ provided that the transfer-matrix associated with  $G_{\tube}(f,g)$ is known. 
 Details of the approach will be published elsewhere.
%are given in Section~\ref{ssec:montecarlo}.  
%
Based on the exact results of Figure \ref{fig_fg},  we focused on the fixed-span model at $g=0$ where the probabilities of the 
smallest trefoil patterns were greatest (compared to the other models). 
Similarly we focus on the $3\times1$ tube, since knots are far more common than in $2\times1$ while the transfer matrices are small enough as to make simulations and enumerations reasonably efficient.

\begin{figure}
\centering
\begin{subfigure}{0.49\textwidth}
\centering
%\resizebox{\textwidth}{!}{\input{3x1_plots/new/plot-3x1-h-3_1-NL&L.tex}}
%\tikzsetnextfilename{figure5left}
%\resizebox{\textwidth}{!}{\input{plot-3_1-LandNL-EL-sHadded.file}}
\resizebox{\textwidth}{!}{\includegraphics{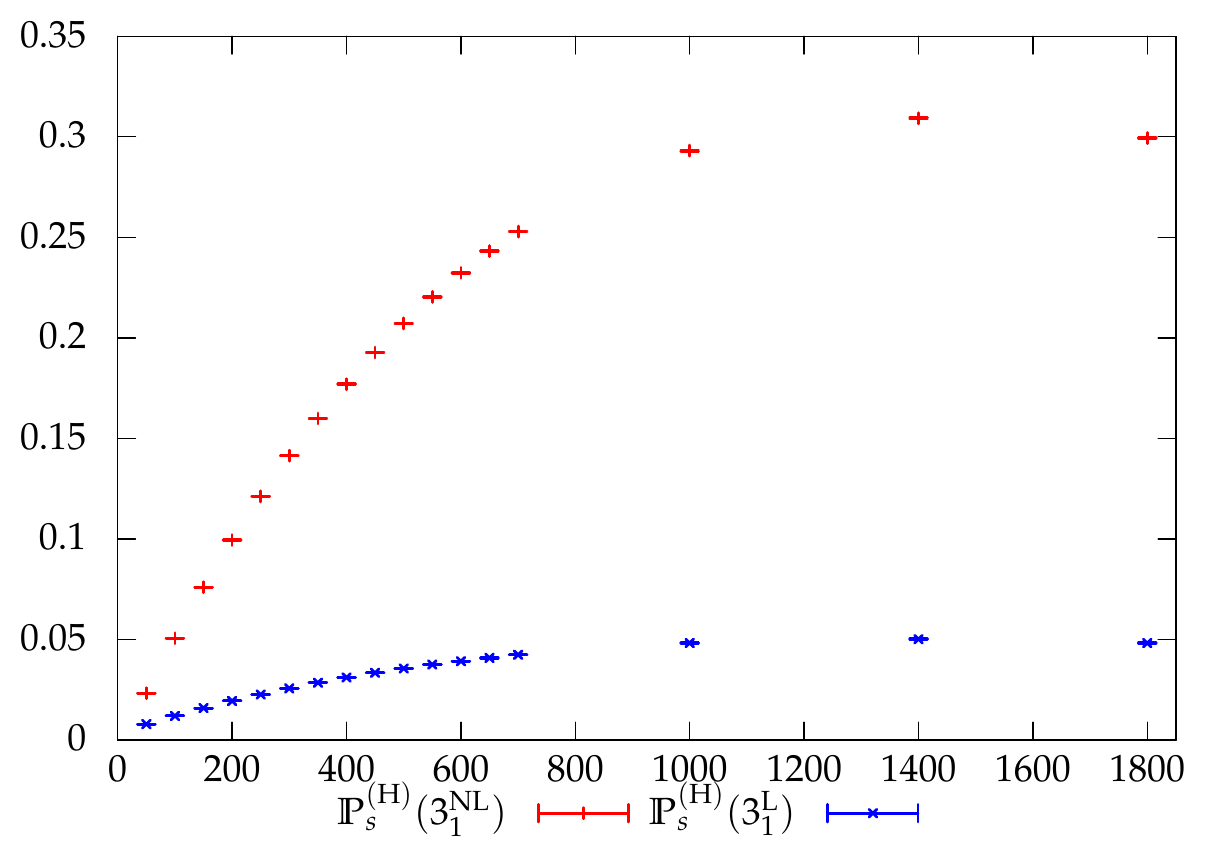}}
\end{subfigure}
\begin{subfigure}{0.49\textwidth}
\centering
%\tikzsetnextfilename{figure5right}
%\resizebox{\textwidth}{!}{\input{plot-3_1-LandNL-EL-sadded.file}}
\resizebox{\textwidth}{!}{\includegraphics{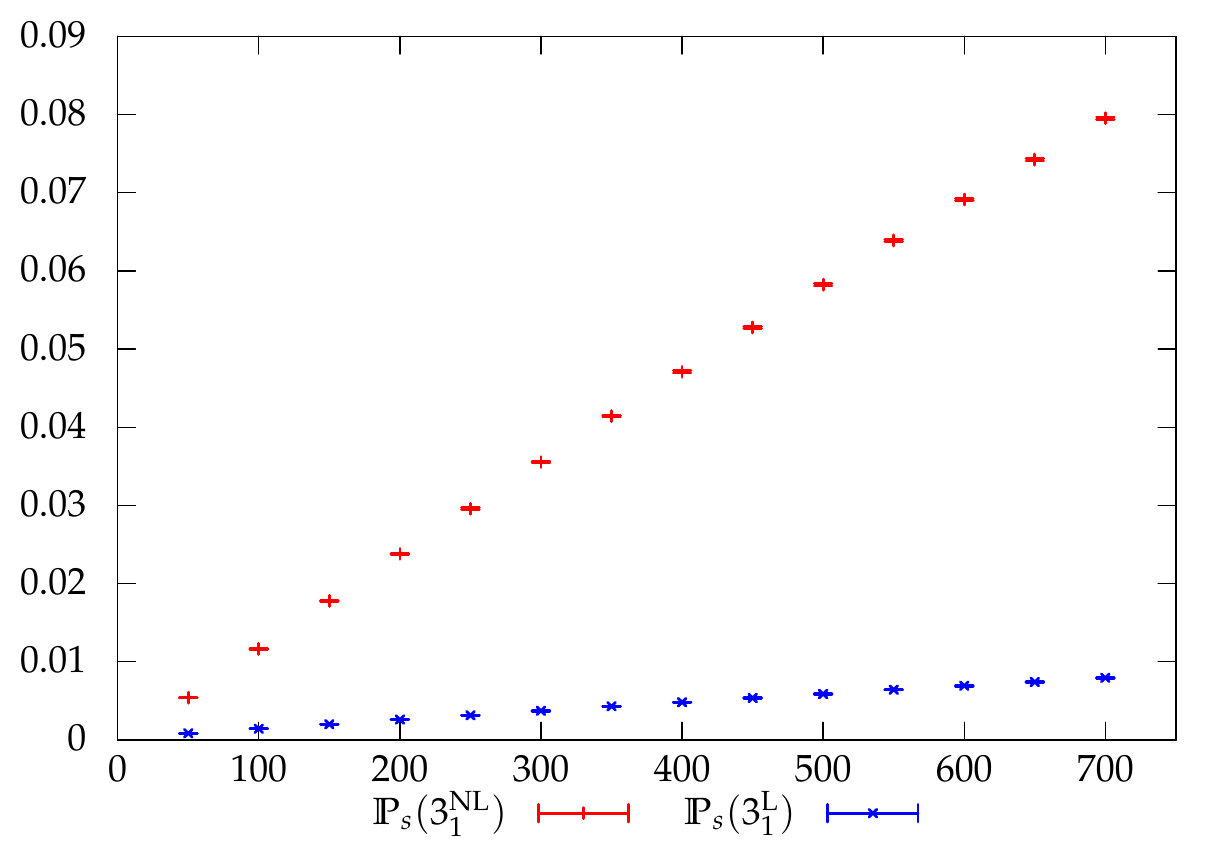}}
\end{subfigure}
\caption{Plots of the probabilities of non-local (red) and local (blue) trefoil knots, for Hamiltonian (above) and all (below) polygons in the $3\times1$ tube, sampled uniformly from the fixed-span ensemble ($g=0$). The horizontal axis is span $s$. The (barely visible) error bars represent 95\% confidence intervals.}
\label{fig:trefoils-local-nonlocal}
\end{figure}

\begin{figure}
\centering
\begin{subfigure}{0.49\textwidth}
\centering
%\resizebox{\textwidth}{!}{\input{3x1_plots/new/plot-3x1-h-4_1-NL&L.tex}}
%\tikzsetnextfilename{figure6left}
%\resizebox{\textwidth}{!}{\input{plot-4_1-LandNL-EL-sHadded.file}}
\resizebox{\textwidth}{!}{\includegraphics{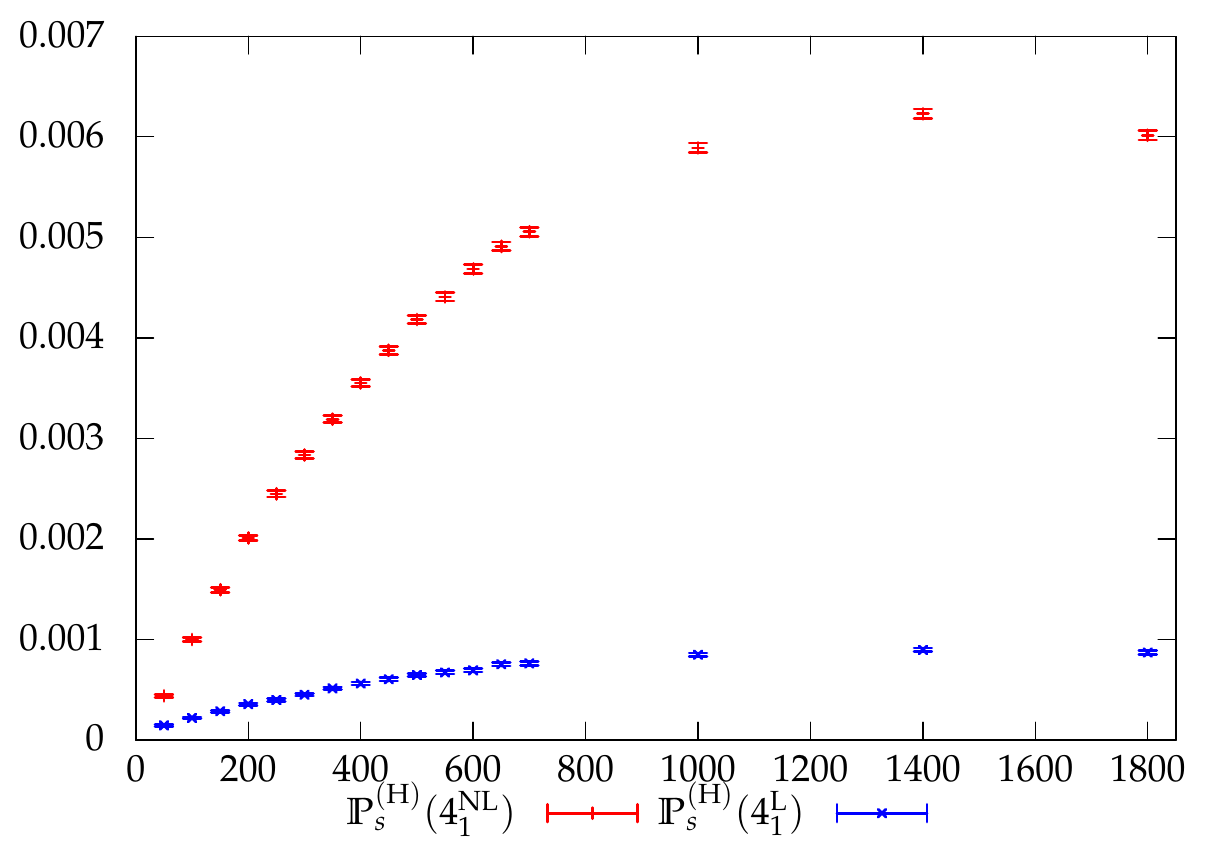}}
\end{subfigure}
\begin{subfigure}{0.49\textwidth}
\centering
%\resizebox{\textwidth}{!}{\input{3x1_plots/new/plot-3x1-a-4_1-NL&L-lines.tex}}
%\tikzsetnextfilename{figure6right}
%\resizebox{\textwidth}{!}{\input{plot-4_1-LandNL-EL-sadded.file}}
\resizebox{\textwidth}{!}{\includegraphics{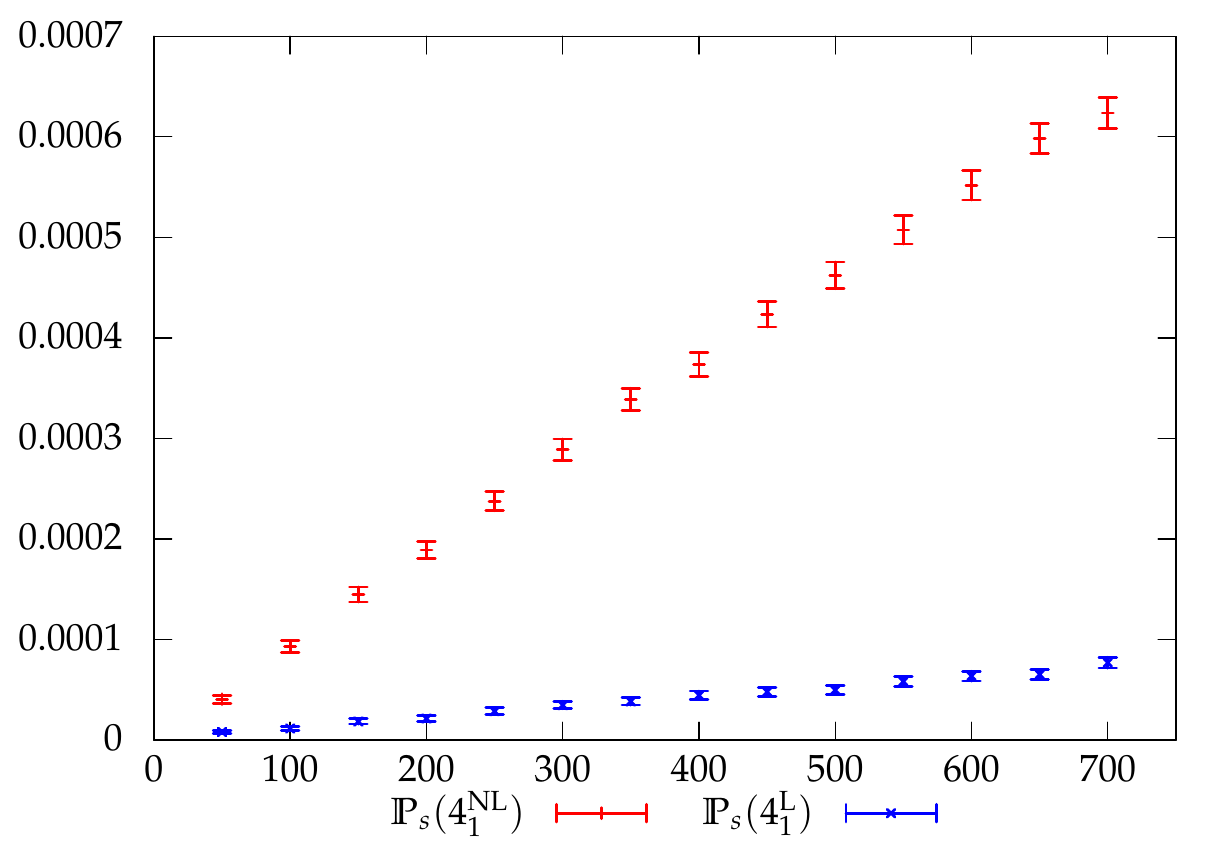}}
\end{subfigure}
\caption{Plots of the probabilities of non-local (red) and local (blue) figure-eight knots, for Hamiltonian (above) and all (below) polygons in the $3\times1$ tube, sampled uniformly from the fixed-span ensemble ($g=0$). The horizontal axis is span $s$. Error bars represent 95\% confidence intervals.}
\label{fig:figureeights-local-nonlocal}
\end{figure}

\begin{figure}
\centering
\begin{subfigure}{0.49\textwidth}
\centering
%\resizebox{\textwidth}{!}{\input{3x1_plots/new/plot-3x1-h-345-NLfrac.tex}}
%\tikzsetnextfilename{figure7left}
%\resizebox{\textwidth}{!}{\input{plot-hratios-EL-Hadded.file}}
\resizebox{\textwidth}{!}{\includegraphics{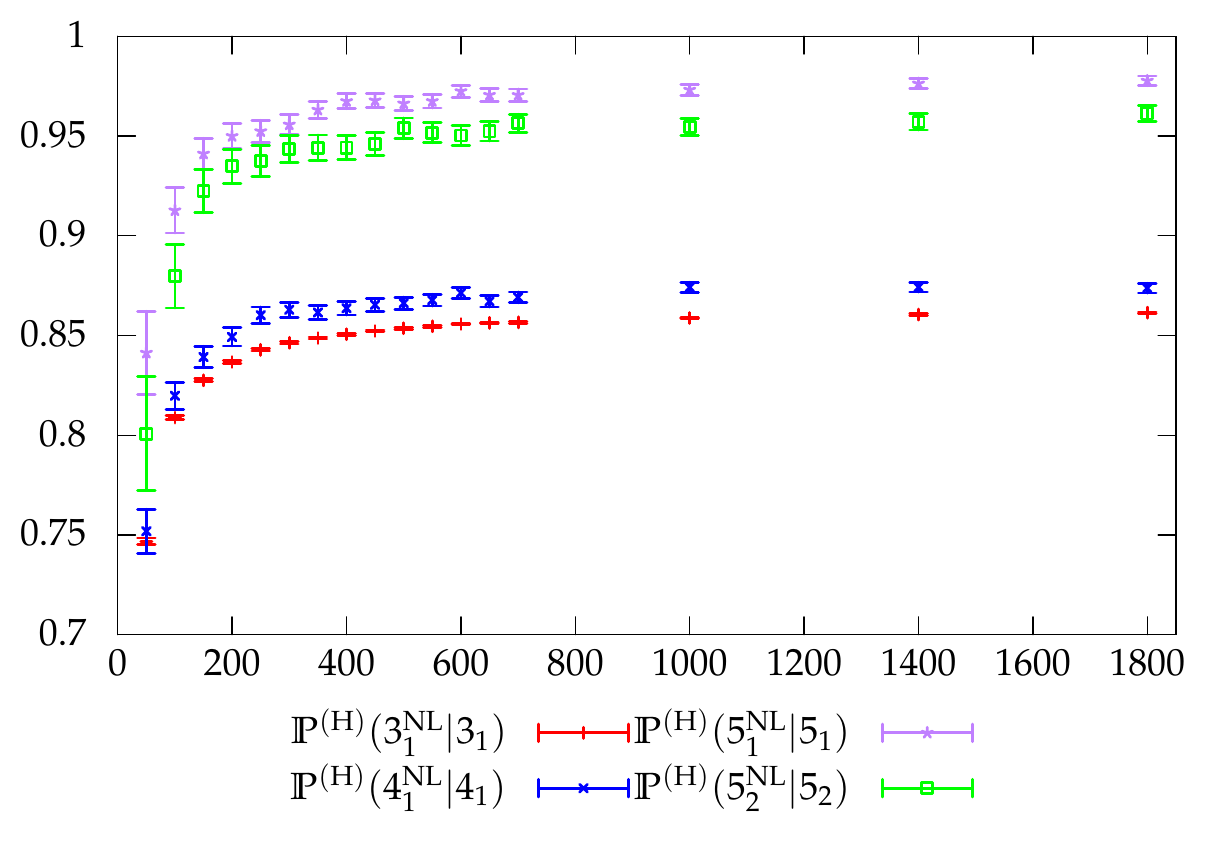}}
\end{subfigure}
\begin{subfigure}{0.49\textwidth}
\centering
%\resizebox{\textwidth}{!}{\input{3x1_plots/new/plot-3x1-a-345-NLfrac-lines.tex}}
%\tikzsetnextfilename{figure7right}
%\resizebox{\textwidth}{!}{\input{plot-aratios-EL.file}}
\resizebox{\textwidth}{!}{\includegraphics{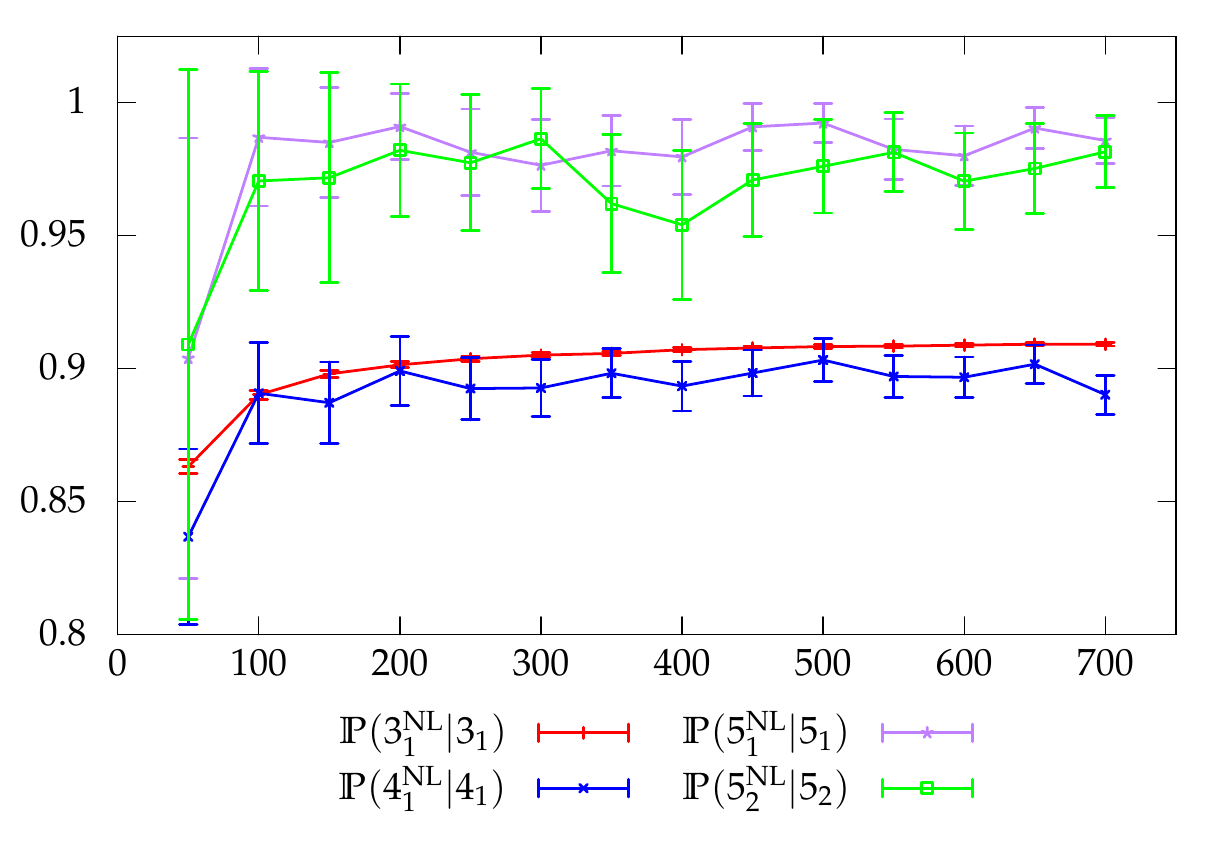}}
\end{subfigure}
\caption{Relative frequencies of non-local knots for trefoils (red), figure-eight knots (blue), $5_1$ knots (purple) and $5_2$ knots (green), for Hamiltonian (above) and all (below) polygons in the $3\times1$ tube, sampled uniformly from the fixed-span ensemble. Note that the lines joining the points in the second plot have been added only to aid the reader, and do not indicate any additional data.}
\label{fig:345-nonlocal-frac}
\end{figure}

For the Monte Carlo results, we begin by investigating the probabilities of some simple knots. In Figure~\ref{fig:trefoils-local-nonlocal} we plot the probabilities of local and non-local trefoils in the fixed-span ensemble, for both Hamiltonian and all polygons. Specifically we plot $\mathbb{P}_s(K^\textup{NL})$ and  $\mathbb{P}_s(K^\textup{L})$ which are respectively the observed proportions of span $s$ polygons which have knot-type $K=3_1$ and contain a non-local (NL) or local (L) knot pattern.  We also plot the corresponding observed proportions of span $s$ Hamilonian polygons: $\mathbb{P}_s^{(\mathrm{H})}(K^\textup{NL}),  \mathbb{P}_s^{(\mathrm{H})}(K^\textup{L})$.
The corresponding data for figure-eight knots ($K=4_1$)  is illustrated in Figure~\ref{fig:figureeights-local-nonlocal}.
The relative frequencies of non-local knots for trefoils ($3_1$), figure-eight knots ($4_1$), $5_1$ and $5_2$ knots are illustrated in Figure~\ref{fig:345-nonlocal-frac}.
For example, the relative frequency of non-local trefoil knots amongst span $s$ trefoil polygons is  $\mathbb{P}_s(3_1^\text{NL}|3_1)=\mathbb{P}_s(3_1^\text{NL})/\mathbb{P}_s(3_1)$.
Although not shown, similar trends were observed for the $2\times 1$ tube. Our observations  lead to the following conclusion.

%Although not shown, similar trends were observed for the $2\times 1$ tube which leads to the following conclusion.

\begin{result}\label{result:MC-probabilities}
For $\tube_{2,1}$ and $\tube_{3,1}$,  based on i.i.d. samples from the limiting distribution of $\mathbb{P}_{s}^{({\rm{sp}}, 0)}$ over a range of spans $s$ (10 million polygons per span),  we observe that the probability of occurrence of a non-local  $3_1$, $4_1$, $5_1$ or $5_2$ knot pattern is significantly greater than that of the corresponding local knot pattern (e.g.~Figures~\ref{fig:trefoils-local-nonlocal} (bottom) and~\ref{fig:figureeights-local-nonlocal} (bottom)). The same holds for Hamiltonian polygons sampled from the limiting distribution of $\mathbb{P}_{s}^{(\mathrm{H})}$ (e.g~Figures~\ref{fig:trefoils-local-nonlocal} (top) and~\ref{fig:figureeights-local-nonlocal} (top)).
Furthermore, for sufficiently long polygons in $\tube_{3,1}$, for each of these knot-types the proportion of non-local patterns amongst all observed knot patterns of that type is greater than 85\% (see Figure
\ref{fig:345-nonlocal-frac}).
\end{result}

\begin{figure}
\centering
\begin{subfigure}{0.49\textwidth}
\centering
%\resizebox{\textwidth}{!}{\input{3x1_plots/new/plot-3x1-hs-3_1-NL&L.tex}}
%\tikzsetnextfilename{figure8left}
%\resizebox{\textwidth}{!}{\input{plot-hs-EL-Hadded.file}}
\resizebox{\textwidth}{!}{\includegraphics{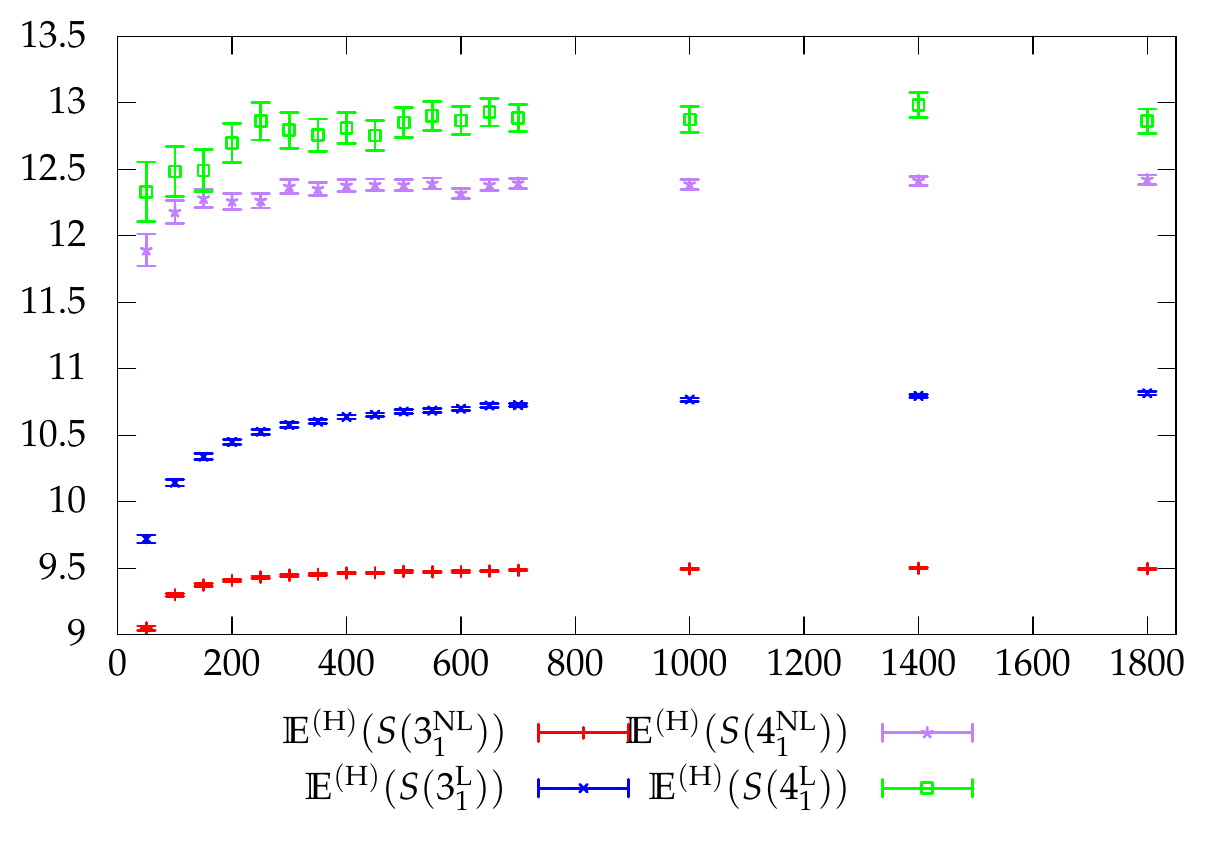}}
\end{subfigure}
\begin{subfigure}{0.49\textwidth}
\centering
%\resizebox{\textwidth}{!}{\input{3x1_plots/new/plot-3x1-as-3_1-NL&L.tex}}
%\tikzsetnextfilename{figure8right}
%\resizebox{\textwidth}{!}{\input{plot-as-EL.file}}
\resizebox{\textwidth}{!}{\includegraphics{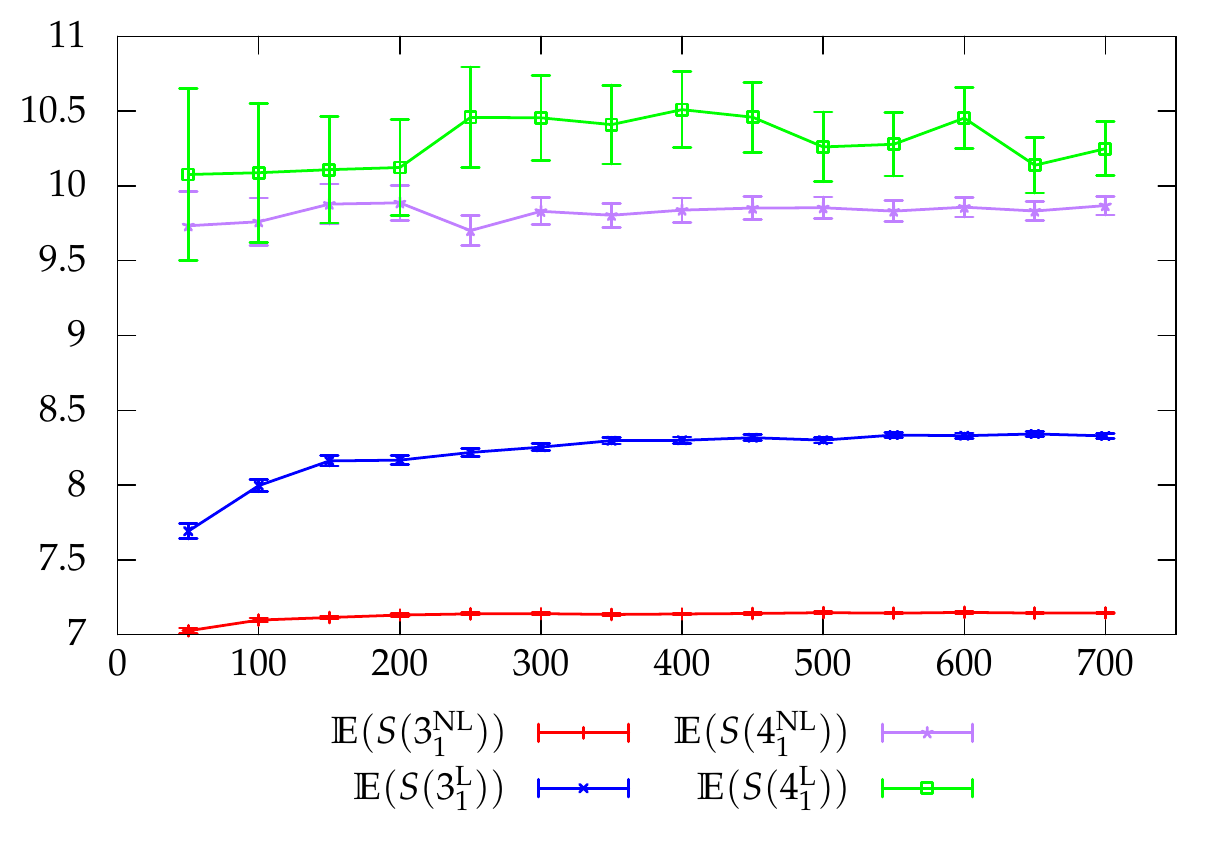}}
\end{subfigure}
\caption{Average span of non-local trefoil (red), local trefoil (blue), non-local figure-eight (purple) and local figure-eight (green)  knot components versus overall polygon span, for Hamiltonian (above) and all (below) polygons in the $3\times1$ tube, sampled uniformly from the fixed-span ensemble.}
\label{fig:3_1-spans}
\end{figure}

Next, we consider the average span of knot patterns within polygons. Recall that we define the span of a knot pattern to include the two 2-sections which bound it on the left and right. In Figure~\ref{fig:3_1-spans} we plot the average span of non-local and local trefoil and figure-eight knot patterns, for both Hamiltonian and all polygons sampled uniformly from the fixed-span ensemble.
Although not shown, similar trends were observed for the $2\times 1$ tube. Our observations lead to the following conclusion.

\begin{result}\label{result:MC-spans}
For $\tube_{2,1}$ and $\tube_{3,1}$,  based on i.i.d. samples from the limiting distribution of $\mathbb{P}_{s}^{({\rm{sp}}, 0)}$ over a range of spans $s$ (10 million polygons per span),  we observe that the average spans of non-local  $3_1$, $4_1$, $5_1$, and $5_2$ knot patterns are smaller than those of the corresponding local knot patterns (e.g.~Figure~\ref{fig:3_1-spans} (bottom)). The same holds for Hamiltonian polygons sampled from the limiting distribution of $\mathbb{P}_{s}^{(\mathrm{H})}$ (e.g.~Figure~\ref{fig:3_1-spans} (top)).
\end{result}

The results of this section provide strong numerical evidence that at least for small tube sizes, non-local knot patterns are more likely and have on average shorter span than local knot patterns.   In the next section we discuss to what extent results related to this can be proved.

\section{Theoretical results}\label{sec:theoretical}

In this section we state results that can be proved related to the occurrence of non-local and local knot patterns in tubes.  We focus on the statement of the results and leave most details of any proofs to the appendix (Section~\ref{ssec:proofs}).  First we present results related to which knot patterns can occur in a tube of a given size and what is known about the minimum span of such knot patterns.
Then we present results on the probability of occurrence of knot patterns.

Because of the strict geometric confinement, not all knots are embeddable in a given $\tube$;  the dimensions of the tube together with the {\em trunk} of the knot (defined in the appendix) determine whether a particular knot is embeddable~\cite{Ishihara_2016}.  
Specifically, it has been shown~\cite{Ishihara_2016} that a knot or link $K$ can be confined to  $\mathbb{T}_{L,M}$
if and only if $\trunk(K)<(L+1)(M+1)$.  This gives that all polygons in a $1 \times 1$ tube are unknots and, in a $2 \times 1$ tube, the knots that can be embedded have been completely classified.  For example, prime knots in the $2 \times 1$ tube are 2-bridge knots.
By applying the method used in~\cite{Ishihara_2016},  we obtain the following results about which knot patterns can occur in $\mathbb{T}_{L,M}$ (see the appendix for some details of the proofs).

\begin{result}\label{result:mxn}
\textbf{(A)}~A knot $K$ admits a proper non-local knot pattern in $\mathbb{T}_{L,M}$ if and only if $\trunk(K)$ $<(L+1)(M+1)$.
\textbf{(B)}~A knot $K$ admits a proper local knot pattern  in $\mathbb{T}_{L,M}$ if $\trunk(K)<(L+1)(M+1)-2$.
\end{result}

Note that this result, while definitive for the case of non-local knot patterns, leaves open the possibility that for a given knot $K$ which is embeddable in $\tubeLM$ there might not be an associated local knot pattern in $\tubeLM$.  For the case of a $2 \times 1$ tube, however, we establish a more definitive result~--~namely, if the knot is embeddable in the tube then there also exist both non-local and local knot patterns for the knot in the tube.  
Figure~\ref{fig:polygon_2x1_hinge_section_block} shows a trefoil polygon in a $2 \times 1$ tube that yields a non-local knot pattern (see Figure~\ref{fig:polygon_2x1_2strings}).  Figure~\ref{fig:polygon_2x1_2stringslocal} shows a local knot pattern in the same tube; note that the span of this local knot pattern is one greater than that shown in Figure~\ref{fig:polygon_2x1_2strings}.
For the $2 \times 1$ tube,  the arguments used in~\cite{Ishihara_2016} can also be extended to prove that this difference in span holds for the smallest knot patterns of knots with up to 5 crossings.  In summary the following result can be proved.

\begin{result}\label{result:2x1_allnonlocal}
Given a prime knot $K\neq 0_1$ that can occur in a $2 \times 1$ tube, there exists at least one proper local knot pattern and at least one proper non-local knot pattern.  Furthermore, at least for $K\in\{3_1,4_1,5_1,5_2\}$,  the span of a smallest proper local knot pattern of $K$ in $\tube_{2,1}$ is greater than that of a smallest proper non-local knot pattern of $K$ in $\tube_{2,1}$.
\end{result}

For the $3\times 1$ tube any knot that can occur in $\tube_{2,1}$ will also have non-local and local knot patterns; see Figure~\ref{fig3_non-localdef} for some examples.
We note that for the $3\times 1$ tube, by an
exhaustive search, we have determined that the span of a smallest local knot pattern of $3_1$ is also one greater than that of a smallest proper non-local knot pattern of $3_1$ (Figures~\ref{fig:local_trefoil} and~\ref{fig:nonlocal_trefoil} are examples of such smallest knot patterns for Hamiltonian polygons in the $3\times 1$ tube). We also looked at the smallest prime knots which cannot occur in $\tube_{2,1}$ but can in $\tube_{3,1}$: $8_5, 8_{10},
8_{15}-8_{21}$.  For all these knots, it is possible to construct a local pattern.  

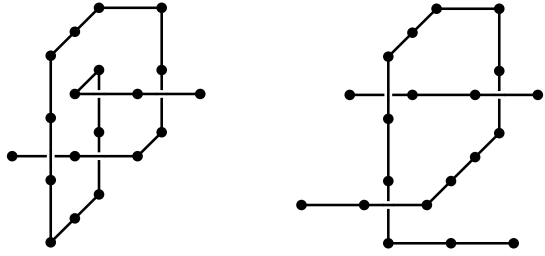
\begin{figure}
\centering
\resizebox{0.45\textwidth}{!}{
\begin{tikzpicture}[rotate around x=270, scale=1.7]
\tikzset{vblack/.style={circle, draw=black, line width=2.5pt, fill=black, inner sep=2pt}}
\begin{knot}[consider self intersections=true, clip width=3, clip radius=4pt, end tolerance=1pt]
\strand[line width=2pt, black] (0,1,1) node [vblack] {} -- (1,1,1) node [vblack] {} -- (2,1,1) node [vblack] {} -- (2,2,1) node [vblack] {} -- (2,2,2) node [vblack] {} -- (2,2,3) node [vblack] {} -- (1,2,3) node [vblack] {} -- (1,1,3) node [vblack] {} -- (1,0,3) node [vblack] {} -- (1,0,2) node [vblack] {} -- (1,0,1) node [vblack] {} -- (1,0,0) node [vblack] {} -- (1,1,0) node [vblack] {} -- (1,2,0) node [vblack] {} -- (1,2,1) node [vblack] {} -- (1,2,2) node [vblack] {} -- (1,1,2) node [vblack] {} -- (2,1,2) node [vblack] {} -- (3,1,2) node [vblack] {};
\flipcrossings{1,3,4}
\end{knot}
\begin{scope}[yshift=-0.4cm]
\begin{knot}[consider self intersections=true, clip width=3, clip radius=4pt, end tolerance=1pt]
\strand[line width=2pt, black] (5,0,1) node [vblack] {} -- (6,0,1) node [vblack] {} -- (7,0,1) node [vblack] {} -- (7,1,1) node [vblack] {} -- (7,2,1) node [vblack] {} -- (7,3,1) node [vblack] {} -- (7,3,2) node [vblack] {} -- (7,3,3) node [vblack] {} -- (6,3,3) node [vblack] {} -- (6,2,3) node [vblack] {} -- (6,1,3) node [vblack] {} -- (6,1,2) node [vblack] {} -- (6,1,1) node [vblack] {} -- (6,1,0) node [vblack] {} -- (7,1,0) node [vblack] {} -- (8,1,0) node [vblack] {};
\strand[line width=2pt, black] (5,2,2) node [vblack] {} -- (6,2,2) node [vblack] {} -- (7,2,2) node [vblack] {} -- (8,2,2) node [vblack] {};
\flipcrossings{2}
\end{knot}
\end{scope}
\end{tikzpicture}
}
\caption{On the left, a shortest knotted arc which can be made a local knot pattern by adding a 3-edge walk. On the right, a shortest non-local knot pattern. Both fit in $\tube_{3,3}$.}
\label{fig:largetref}
\end{figure}

The fact that the span of smallest non-local knot patterns is smaller than that of local knot patterns provides a partial explanation for Result~\ref{result:MC-spans}, that the average span of a non-local knot pattern is smaller than that of a corresponding local knot pattern for small tube sizes.  However, for a large enough tube size, the span of smallest non-local and local knot patterns are expected to be the same; see Figure~\ref{fig:largetref} which shows, on the left, a shortest trefoil arc which can be turned into a span-3 local trefoil knot pattern by adding a 3-edge walk within $\tube_{3,3}$, and on the right is a smallest non-local trefoil knot pattern in $\tube_{3,3}$ also having span 3.  Note, however, for this case the number of edges in the local knot pattern is greater than that for the non-local pattern; this may lead to non-local trefoil patterns being more likely to occur than  local ones even in larger tube sizes, i.e. that Result~\ref{result:MC-probabilities} could continue to hold.  

%\begin{figure}
%\mbox{{\includegraphics[width=0.10\textwidth]{tighttref2}}. $~~$ {\includegraphics[width=0.24\textwidth,angle=0]{3_1nonlocaldrawing}}}
%\caption{On the left, a shortest knotted arc which can be made a local knot pattern by adding a 3-edge walk. On the right, a shortest non-local knot pattern. Both fit in $\tube_{3,3}$.}
%\label{fig:largetref}
%\end{figure}

We next  
discuss some results about the likelihood of
each type of pattern.

There are known ``pattern theorems'' available for both the fixed-edge and fixed-span models studied here (see~\cite{Atapour2008PhD,Atapour_2009}), as well as for Hamiltonian polygons (see~\cite{Eng2014MSc}).   
The theorems focus on proper polygon patterns (see~\cite{Beaton_2016} for more precise definitions) which include the proper knot patterns defined here. 
Given a model and a proper pattern $P$ which can occur in a polygon of the model in $\tube$,   a pattern theorem establishes that there exists an $\epsilon_P>0$  such that, for $n$ sufficiently large, all but exponentially few $n$-edge polygons contain more than $\epsilon_Pn$ copies of $P$.
From such theorems it is known that the knot-complexity of polygons grows as polygon ``size'' grows (size could be measured in terms of edges or span), so that a typical polygon will have a highly-composite knot-type $K=K_1\# K_2 \# \ldots\# K_r$~\cite{SSW}.  Different prime components of the knot could occur as knot patterns in the polygon in a variety of ways.  Our interest here is to investigate how often they are occurring as ``local'' knots versus non-locally.    

To define ``local'' knotting requires the definition of a knot-size measure 
(we have given two possible measures: arclength  and connect-sum knot-size) 
but also a comparison of knot-size to polygon size. 
Here we consider that a polygon's size $m$ is growing without bound and say that it is non-locally knotted with respect to arclength knot-size if at least one of the knots $K_i$ in its prime knot decomposition  has arclength knot-size $a_{K_i} = O(m)$.  In this case we say $K_i$ occurs non-locally in the polygon or is non-locally knotted, and otherwise we say that $K_i$ occurs locally or is locally knotted.  Thus a polygon can be both non-locally and locally knotted (with respect to arclength knot-size) depending on the occurrence of each of its prime components.   
Corresponding definitions can apply to the case of the connect-sum knot-size.  However, to distinguish this case, we say a polygon is ``loosely'' knotted (or contains a loose knot) if at least one of the knots $K_i$ in its prime knot decomposition has connect-sum size $b_{K_i}=O(m)$; otherwise, we say $K_i$ occurs as a tight knot or tightly.  

To explain these definitions more clearly, consider the non-local trefoil pattern $P$ of Figure~\ref{fig:polygon_2x1_2strings}.  If this pattern occurs
in a trefoil polygon $\pi$, then $\pi$ must be formed from a connect sum of unknotted polygons with $P$.  Suppose the size of $\pi \gg$ the span of $P$. If $P$ occurs near the left end (or right end) of $\pi$,  then the arclength knot-size will be short compared to the size of $\pi$ and even though we have classified $P$ as a non-local knot pattern,  we will say that
$P$ has occurred in a ``local" way in $\pi$ and that $\pi$ is locally knotted.   If instead $P$ occurs in the ``middle" of $\pi$ (i.e.~half-way along the span) then we will say that $P$ has occurred in a ``non-local" way in $\pi$ and that $\pi$ is non-locally knotted (since the arclength size of the knot is proportional to the size of $\pi$).   In contrast,  because the connect-sum knot-size of $P$ is small compared to the size of $\pi$,  no matter where it occurs in $\pi$ we will say that $3_1$ occurs as a tight knot.   
On the other hand,  for the local trefoil pattern of Figure~\ref{fig:polygon_2x1_2stringslocal},  no matter where it occurs in a very large sized trefoil polygon $\pi$,  it will always have 
a small arclength knot-size as well as a small connect-sum knot-size and hence $\pi$ will be considered to be both locally and tightly knotted.  

Since the pattern theorems hold for proper knot patterns, they   
tell us that for each of the models in question, any proper knot pattern which can occur in a polygon (or Hamiltonian polygon) will occur with a positive density as polygon size grows (where, depending on the model,  size is measured by edges or by span).  
In particular for any $\tube$,  polygons which contain no knot patterns 
are exponentially rare. Furthermore we have, by considering the local trefoil knot pattern of Figure~\ref{fig:polygon_2x1_2stringslocal} (or an inflated version for the Hamiltonian polygon cases),
%(CHANGE THESE TO REFER TO 2x1 EXAMPLES),
 the following result.

\begin{result}\label{result:all_but_exp_few_locally_knotted}
All but exponentially few sufficiently long polygons in $\tube$ with $M\geq L\geq 2$ or $M\geq 2, L=1$  for any of the three models above are both locally and tightly knotted.
\end{result}

Note that this does not preclude the same polygons from being non-locally or loosely knotted  -- we only know that the knot-types of the polygons are highly complex and that some of the knots  in the knot decomposition will be local trefoil knot patterns as in Figure~\ref{fig:polygon_2x1_2stringslocal}.

Indeed we use the same pattern theorems next to establish there will be non-local knot patterns that occur in a non-local way in all but exponentially few polygons.    Let $P$ be the non-local trefoil pattern in Figure~\ref{fig:polygon_2x1_2strings} (or expanded appropriately to fill the tube in the Hamiltonian polygon case) and let $\epsilon_P>0$  
be as needed for the appropriate pattern theorem.  Let $n$ be sufficiently large that $\epsilon_Pn > 4$.  We consider any $n$-edge polygon that contains at least $\epsilon_P n$ translates of $P$.  Since the span of $P$ is 6, and since $P$ cannot overlap itself,  the $(\epsilon_P n/2)$th copy of $P$ must start at a constant $x$ plane with $x\geq (6\epsilon_Pn/2)-1$   so that the arclength knot-size of $3_1$ in the polygon must be greater than or equal to
  $6\epsilon_Pn=O(n)$.  Thus $3_1$ occurs non-locally at least once in such a polygon.   By the pattern theorem, all other polygons are exponentially rare.   
  Thus we have the following result.

\begin{result}\label{result:all_but_exp_few_nonlocally_knotted}
All but exponentially few sufficiently  long polygons in $\tube$ with $M\geq L\geq 2$  or $M\geq 2, L=1$ are non-locally knotted.
\end{result}
 
Meanwhile for self-avoiding walks in $\tube$,  at least for $f=0$, the scenario depicted in Figure~\ref{fig:nonlocal_loop_walk}, in which a non-local  knot pattern occurs in a walk in a non-local way, is exponentially rare.  To see this, first note that any polygon can be turned into a walk by removing one edge.  Thus from~\eqref{eqnmudef} for $f=0$, we know that polygons are exponentially rare in the set of walks.   Next consider the subset of $n$-step walks in $\tube$ which contain a non-local knot-type $K$ pattern at a location in the walk such that the arclength of  $K$ is $\alpha n$ for some $\alpha >0$ (as in Figure~\ref{fig:nonlocal_loop_walk}).   Each such walk can be  decomposed into a polygon with at least $\alpha n$ edges (and having $K$ in its knot-decomposition) and a walk with length at most $(1-\alpha n)$.
Thus this subset of walks will have an exponential growth rate which is strictly less than all walks in $\tube$.  This gives the following result.
 
\begin{result}\label{result:SAWs-nonlocal-exp-rare}
Self-avoiding walks which contain a non-local knot pattern in a non-local way are exponentially rare in the set of all walks in $\tube$.
\end{result}

(Note that this result does not contradict the pattern theorem for walks in $\tube$ proved in~\cite{Soteros_1988} because the proper knot patterns defined here are not examples of proper walk patterns.)

\section{Conclusion}\label{sec:conclusion}

We have used self-avoiding polygons to model ring polymers confined to narrow tubes. For this model we have used a standard approach for measuring the size of a knot to define  ``local'' and ``non-local'' knotting.

We have then provided both theoretical and numerical evidence that when ring polymers are confined to very narrow tubes, at equilibrium and assuming all states are accessible, non-local knotting is more likely than local knotting. One reason for this may be that non-local knot configurations, at least for the simplest knots, are on average smaller than their local counterparts. In small tube sizes they are smaller both in span and in edge count.  

These results can be compared and contrasted with recent numerical models~\cite{Suma_2017} for DNA knots in a solid-state nanopore~\cite{Plesa_2016}. Two different modes of knot translocation have been found, which closely correspond with our definitions of local and non-local knotting (see Figures 4(C) and (D) in~\cite{Suma_2017}). For that (non-equilibrium) model, the relative frequencies of these modes (measured as the knot passes through the pore) depended on the length of the DNA molecule and the initial conditions prior to translocation.

We also provided theoretical evidence that for linear chains, non-local knotting is exponentially rare, due to the entropic disadvantage of a long bend. This is comparable to what has been observed experimentally~\cite{Plesa_2016}.

\section{Appendix}\label{sec:appendix}

\subsection{Comments on characterisation of local knot patterns}\label{ssec:comments}
\begin{figure}
\centering
\includegraphics[width=0.6\textwidth]{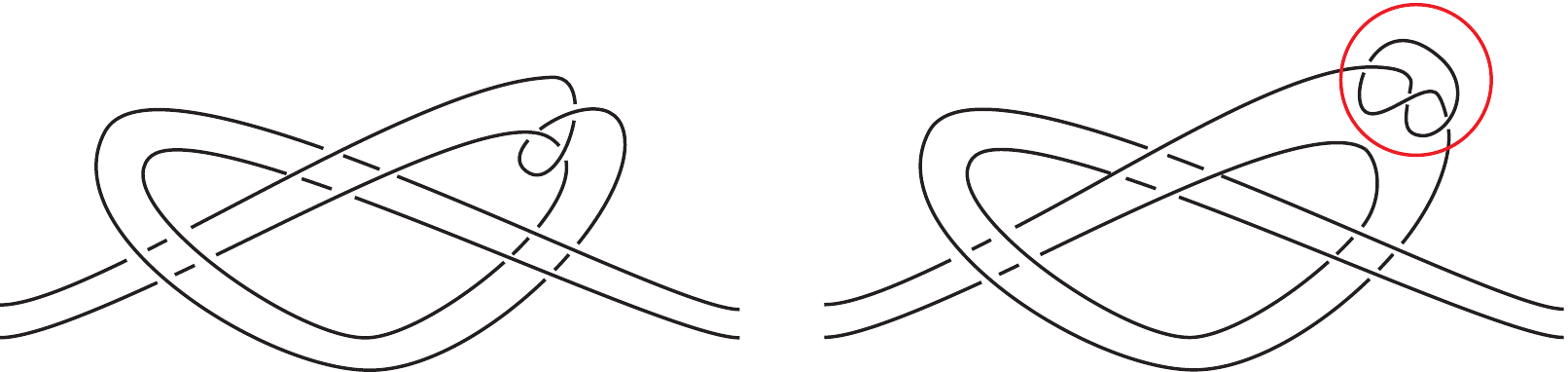}
\caption{Two examples of local knot patterns. The right one has a red 2-sphere intersecting the pattern in two points and surrounding a local knot.}
\label{fig:versus-2}
\end{figure}

In this paper we have defined a classification scheme for local knot patterns that is relatively simple to implement for polygons in tubes, however, this classification is not necessarily consistent with the topological definition of local knot.  In particular, a topologically precise approach to defining a local knot pattern is as follows.  
For a given knot pattern,  if we can find a 2-sphere that intersects the pattern in exactly two points and such that it surrounds a knotted arc, then the knot pattern contains a local knot.  Figure~\ref{fig:versus-2} right shows an example of a knot pattern with such a 2-sphere.   However, by our classification scheme, 
Figure~\ref{fig:versus-2} left shows an example of a knot pattern $\sigma$ which we are currently classifying as a local knot pattern even though there is no 2-sphere satisfying the property just described.   In particular, $DC(\sigma)=NC_1(\sigma)=NC_2(\sigma)=3_1$, however, there is no 2-sphere that isolates (by surrounding it at the exclusion of anything else) either of these trefoils and intersects the pattern in only two points.  

Because the number of edges needed to create patterns such as that shown in Figure~\ref{fig:versus-2} (left) is large, such patterns are not very likely to occur and we do not believe that their existence will affect our conclusions.

\subsection{Outlines of Proofs of Theoretical Results}\label{ssec:proofs}

Given a polygon $\pi\in\mathcal{P}_\tube$, a \emph{hinge} $H_k$ of $\pi$ is the set of edges and vertices lying in the intersection of $\pi$ and the $y$-$z$ plane defined by $\{(x,y,z):x=k\}$. See Figure~\ref{fig:polygon_2x1_hinge_section_block} for an example.

\begin{figure}
\centering
\begin{subfigure}{0.48\textwidth}
\centering
\includegraphics[scale=0.4]{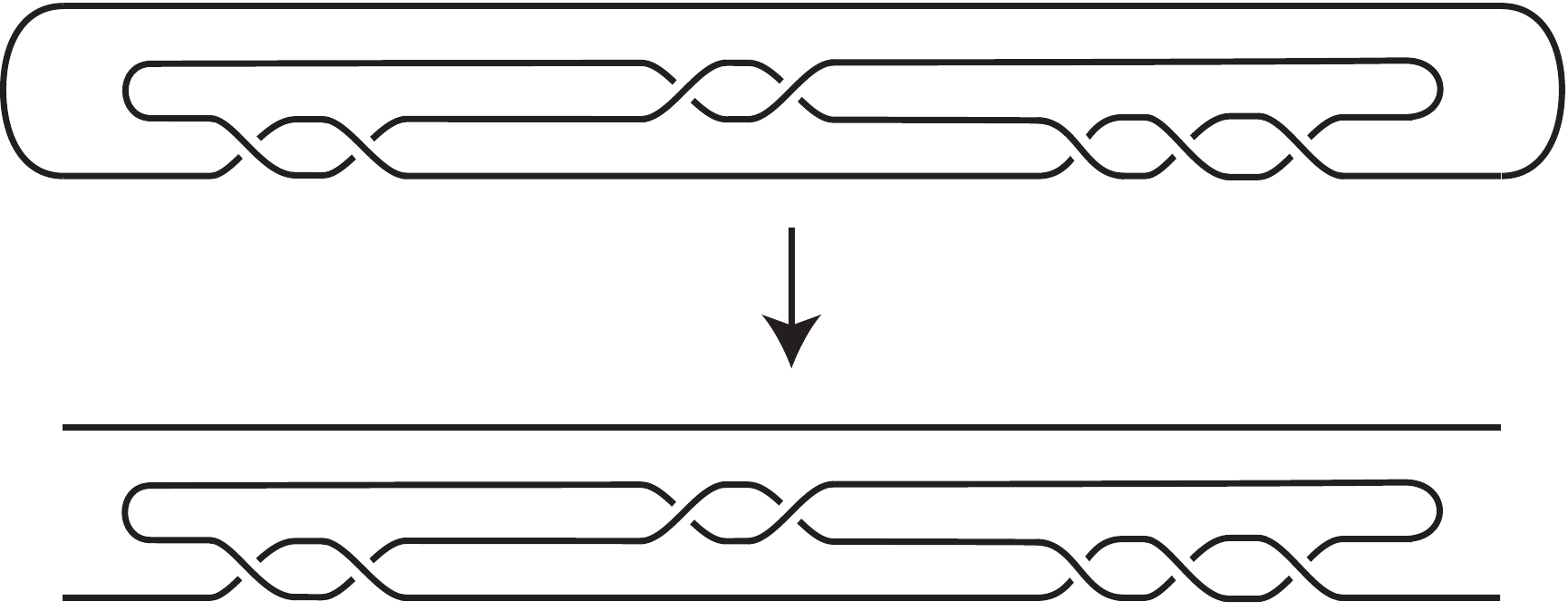}
\caption{}\label{fig:local2bridge}
\end{subfigure}
\begin{subfigure}{0.48\textwidth}
\centering
\includegraphics[scale=0.5]{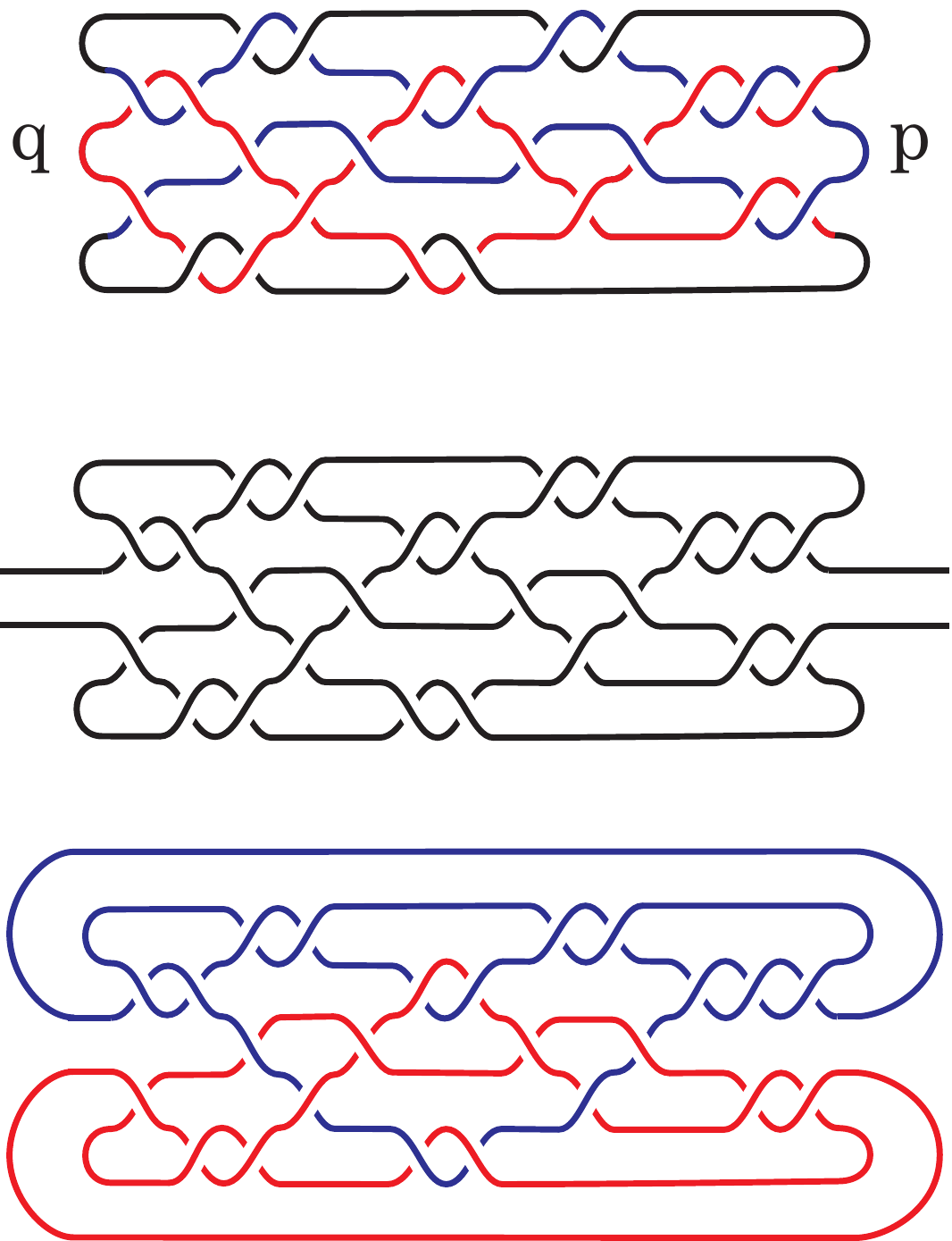}
\caption{}\label{fig:trunk6}
\end{subfigure}
\begin{subfigure}{0.6\textwidth}
\centering
\includegraphics[width=\textwidth]{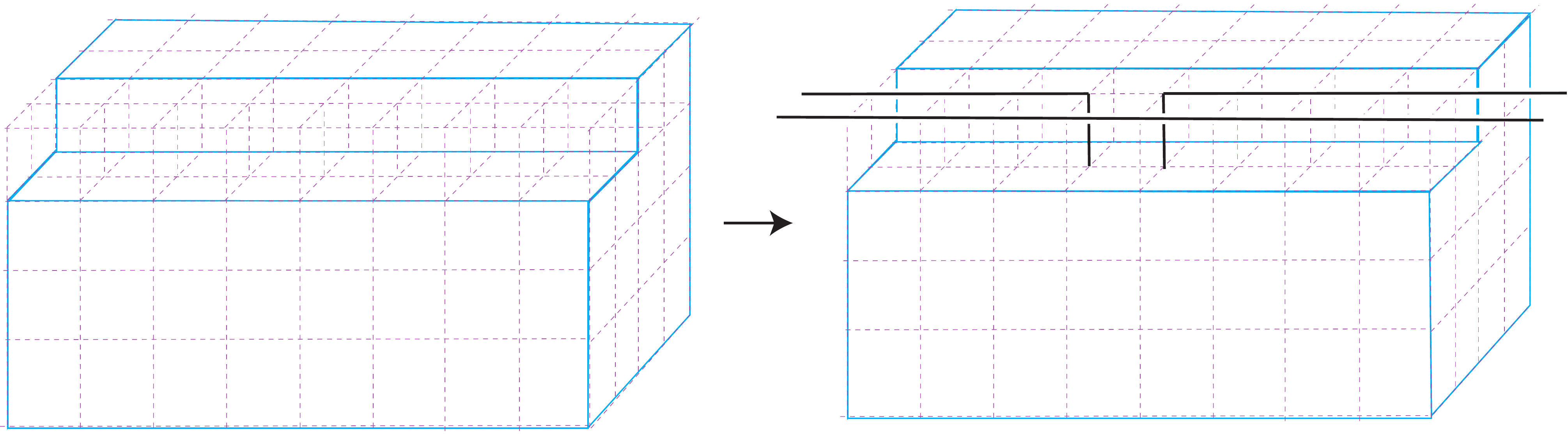}
\caption{}\label{fig:localknotgeneral}
\end{subfigure}
\caption{
\subref{fig:local2bridge} A local $7_5$ pattern obtained by opening ends of Conway's normal form.
\subref{fig:trunk6} When $\trunk(K)\ge 6$ we can choose $p$ and $q$ so that each arc of $\pi_K-\{p,q\}$ contains at least two critical points.
%; 5$7_5=C(2,2,3)$ is arranged to $C(3,-1,-2,1,2)$ which can be constructed in $2\times 1$ tube.
\subref{fig:localknotgeneral} By pulling out a part, we can construct a local knot pattern.
}\label{fig:result1}
\end{figure}

\begin{proof}[Proof of Result \ref{result:mxn}]
\textbf{(A)} The {\em trunk} of a knot or link $K$ is an invariant defined by
$\trunk(K)=\min_{E}\max_{t\in \mathbb R} |h^{-1}(t)\cap E|$,
where $E$ is any embedding of $K$ in $\mathbb R^3$ and $h:\mathbb R^3\to \mathbb R$ is any given height function~\cite{Ozawa2010,Ishihara_2016}. For another invariant, the bridge number $b(K)$ of $K$, $\
trunk(K)$ satisfies $\trunk(K) \le 2b(K)$.

By~\cite[Theorem 1]{Ishihara_2016}, we can construct a polygon of knot type $K$ in $\mathbb{T}_{L,M}$ if and only if $\trunk(K)<(L+1)(M+1)$.
Then we can obtain a proper knot pattern from such a polygon by opening its ends, i.e. by removing an edge or edges (as appropriate) in each of the left-most and right-most hinges.
See Figure~\ref{fig:local2bridge}.
We will show that there is a polygon which can be opened at each end to yield a proper non-local knot pattern.

First we consider the case where $\trunk(K)\ge 6$.
Take a height function $h$ and an embedding of $K$, $\pi_K$, in $\tube$ such that  $\trunk(K)$ is attained and such that $\pi_K$ has the minimal number of critical points with respect to $h$.
We can choose one maximal point $p$ and one minimal point $q$ to make a proper knot pattern so that each of the two arcs of $\pi_K-\{p,q\}$ has at least two critical points.
See Figure~\ref{fig:trunk6}.
Let $K_1$ and $K_2$ be the components of the link obtained by taking the numerator closure of $\pi_K-\{p,q\}$.
Then neither of $K_1$ nor $K_2$ is $K$ by the minimality of the number of critical points of $\pi_K$.
It follows that  the pattern is non-local.
We can construct a polygonal model of $K$ satisfying the above conditions in a given tube.

Suppose $\trunk(K)=4$.
First we consider the case where $K$ is prime, i.e., $K$ is a 2-bridge knot.
Take a Conway's normal form with the minimal crossing number.
Then there are at least two strings of the 4-braid corresponding to the Conway's normal form that contain crossings.
%WHAT ARE $K_1$ AND $K_2$?
Then we can make a proper knot pattern so that  both $K_1$ and $K_2$, the components of the numerator closure, contain one each of such strings.
Then the crossing numbers of $K_1$ and $K_2$ are strictly less than that of $K$.
Hence neither of $K_1$ nor $K_2$ is $K$.
We can construct a polygonal model of $K$ satisfying the above conditions in a given tube and it gives a non-local pattern.
Suppose $K$ is a composite knot.
Let $L_1$ and $L_2$ be knots such that $K=L_1\# L_2$ and $L_1$ is a prime knot.
Then by the above argument, we can create a non-local pattern for $L_1$.
By a connected sum operation, we can then construct a polygon of $K$ that gives a non-local proper knot pattern.

\textbf{(B)} Suppose $\trunk(K$) $<(L+1)(M+1)-2$.
Then by using a method of~\cite[Theorem 1]{Ishihara_2016}, we can construct a polygon inside a region in $\mathbb{T}_{L,M}$ as in Figure~\ref{fig:localknotgeneral} (left).
Then by pulling out a part as in Figure~\ref{fig:localknotgeneral} (right) we have a local proper knot pattern. 
\end{proof}

\begin{figure}
\centering
\begin{subfigure}{0.48\textwidth}
\centering
\includegraphics[scale=0.35]{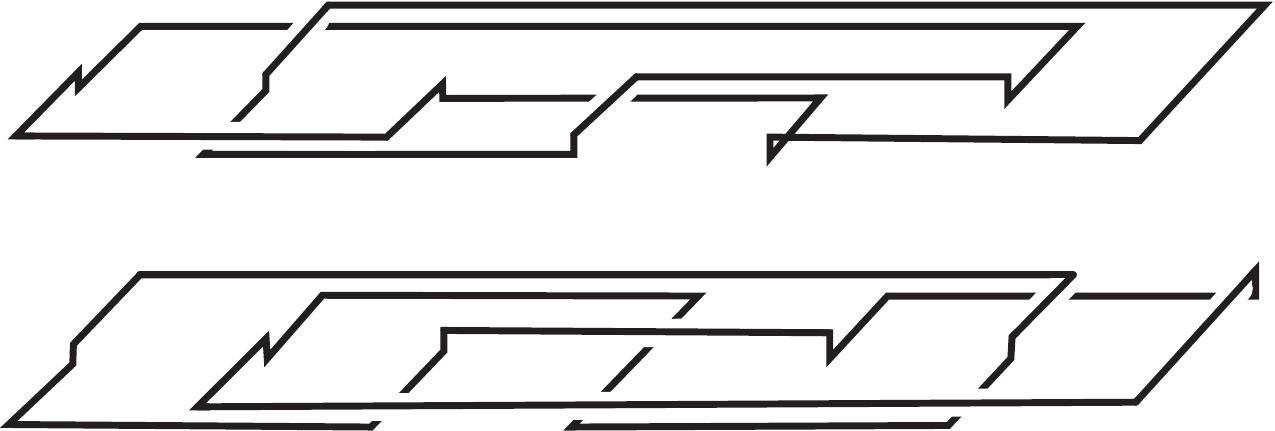}
\caption{}\label{fig:local3_1}
\end{subfigure}
\begin{subfigure}{0.48\textwidth}
\centering
\includegraphics[scale=0.35]{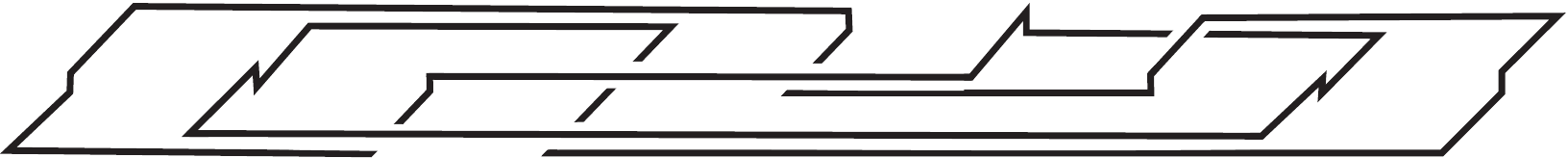}
\caption{}\label{fig:local4_1}
\end{subfigure}
\begin{subfigure}{0.48\textwidth}
\centering
\includegraphics[scale=0.35]{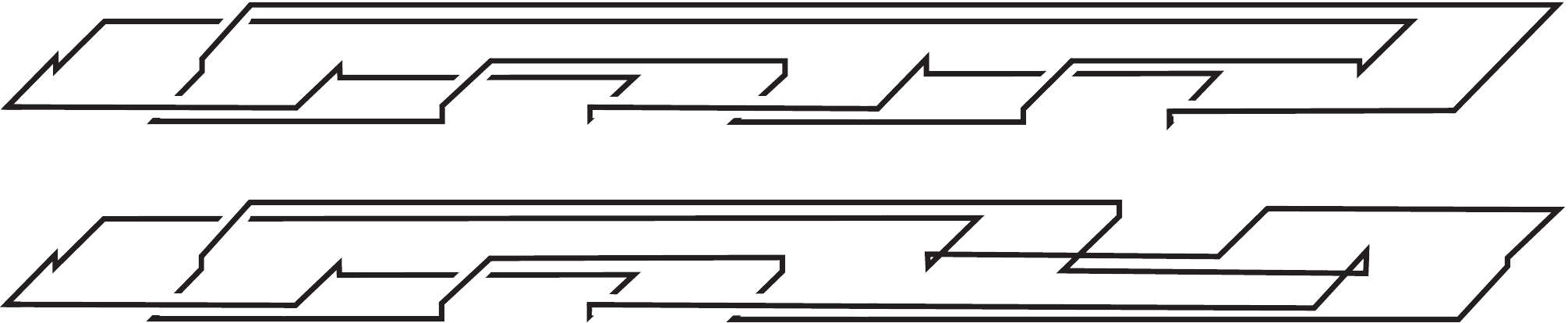}
\caption{}\label{fig:local5_1}
\end{subfigure}
\begin{subfigure}{0.48\textwidth}
\centering
\includegraphics[scale=0.35]{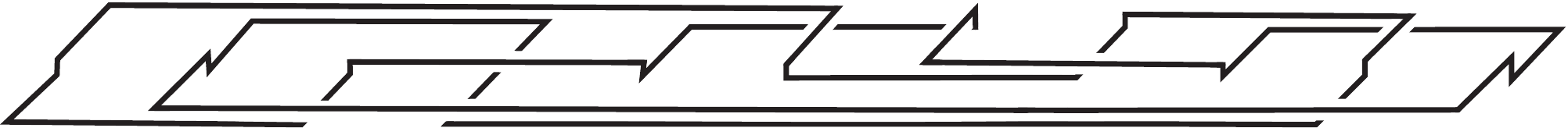}
\caption{}\label{fig:local5_2}
\end{subfigure}
\caption{
\subref{fig:local3_1} Two polygons of $3_1$ in $2\times 1$ tube with the smallest span $6$; the first consists of $36$ edges and the second consists of $38$ edges.
\subref{fig:local4_1} A polygon of $4_1$ in $2\times 1$ tube with the smallest span $8$; this consists of $50$ edges. 
\subref{fig:local5_1}Two polygons of $5_1$ in $2\times 1$ tube with the smallest span $10$; these consist of $60$ edges.
\subref{fig:local5_2} A polygon of $5_2$ in $2\times 1$ tube  with the smallest span $10$; this consists of $62$ edges. 
}\label{fig:result2}
\end{figure}

\begin{proof}[Proof of Result \ref{result:2x1_allnonlocal}]
Any prime knot that can occur in a $2\times1$ tube is $2$-bridge~\cite{Ishihara_2016}.
It is well known that any $2$-bridge knot is represented by Conway's normal form $C(a_1,\cdots,a_n)$, 
which is a closure of a $4$-braid using only two generators $\sigma_1$ and $\sigma_2$~\cite{Murasugi_1996}. 
Since there is no $\sigma_3$ and the fourth string in the Conway's normal form is straight, we have a local knot pattern by opening both ends as in Figure~\ref{fig:local2bridge}.
%Moreover, we can naturally construct such a local proper knot pattern in $2\times 1$ tube, by arranging the Conway's normal form so that all even-numbered integers $a_2,a_4,\ldots$ become $\pm1$, see Fig. \ref{fig:local2bridge} for the case of the $7_5$-knot.
%From a local proper knot pattern, we also have a non-local proper knot pattern for the same knot type by adding twists of two strings on one of the ends as in Fig. \ref{}.

By~\cite[Lemma 3(1)]{Ishihara_2016},
from a knotted polygon in  a $2\times 1$ tube with the smallest span we can obtain a proper knot pattern with the smallest span in the $2\times 1$ tube for that knot type 
%we have a knot polygon in $2\times 1$ tube with the smallest span 
by opening both ends of the polygon. 
For $K\in\{3_1,4_1,5_1,5_2\}$, by applying the argument of~\cite[Theorem 4]{Ishihara_2016}, we can completely characterise the configurations of $K$ with smallest span,
see Figures~\ref{fig:result2}\subref{fig:local3_1}, \subref{fig:local4_1}, \subref{fig:local5_1}, \subref{fig:local5_2} for examples.
We can then conclude that the resulting proper knot patterns are all non-local.
On the other hand, in these cases a local proper knot pattern can be constructed by increasing the span by one by the same method as in Figure~\ref{fig:polygon_2x1_2stringslocal}.
\end{proof}

\section*{Acknowledgments}

CES acknowledges support in the form of a Discovery Grant from NSERC (Canada) and a CPU allocation from Compute Canada's WestGrid.  KI is partially supported by JSPS KAKENHI Grant Number 17K14190. KS is partially supported by JSPS KAKENHI Grant Numbers 26310206, 16H03928, 16K13751 and 17H06463. NRB received support from the PIMS Collaborative Research Group in Applied Combinatorics, and from Australian Research Council grant DE170100186. CES acknowledges helpful discussions with K.~Millett. The authors also acknowledge assistance from Rob Scharein with KnotPlot and that some figures were produced using Rob Scharein's KnotPlot.

\bibliographystyle{plain}
\bibliography{non-local-plain}

%\printbibliography

\end{document}